\documentclass[11pt,a4paper,leqno]{amsart}
\usepackage{fullpage,amsmath,amsthm,amsfonts,amscd,
latexsym,amssymb,eucal,curves}
  \usepackage[all]{xy}
     \usepackage[dvips]{graphicx}
\newcommand{\ZZ}{\mathbb{Z}}

\newcommand{\CC}{\mathbb{C}}

\newcommand{\PP}{\mathbb{P}}
\newcommand{\NN}{\mathbb{N}}

\newcommand{\HH}{\mathbb{H}}

\newcommand{\LL}{\mathbb{L}}

\newcommand{\MM}{\mathbb{M}}

\newcommand{\WW}{\mathbb{W}}

\newcommand{\QQ}{\mathbb{Q}}
\newcommand{\RR}{\mathbb{R}}

\newcommand{\FF}{\mathbb{F}}
\newcommand{\VV}{\mathbb{V}}

\newcommand{\Aut}{\mathop{\rm Aut}\nolimits}

\newcommand{\id}{{\rm id}}
\renewcommand{\Re}{\mathop{\rm Re}}
\renewcommand{\Im}{\mathop{\rm Im}}
\newcommand{\dd}{{\rm d}}

\newcommand{\Gal}{{\rm Gal}}

\newcommand{\cLL}{{\mathcal L}}
\newcommand{\cOO}{{\mathcal O}}
\newcommand{\cEE}{{\mathcal E}}
\newcommand{\cFF}{{\mathcal F}}

\newcommand{\cXX}{{\mathcal X}}
\newcommand{\cYY}{{\mathcal Y}}
\newcommand{\cZZ}{{\mathcal Z}}
\newcommand{\cPP}{{\mathcal P}}
\newcommand{\cUU}{{\mathcal U}}

\newcommand{\SL}{{\rm SL}}

\newcommand{\diag}{{\rm diag}}

\newcommand{\tr}{{\rm tr}}
\newcommand{\sms}{\smallsetminus}

\newcommand{\Su}{S_{\scriptstyle \rm u}}

\newcommand{\ol}{\overline}

\newcommand{\PSL}{{\rm PSL}}

\newtheorem{Defi}{Definition}[section]

\newtheorem{Rem}[Defi]{Remark}

\newtheorem{Prop}[Defi]{Proposition}
\newtheorem{Lemma}[Defi]{Lemma}
\newtheorem{Cor}[Defi]{Corollary}
\newtheorem{Thm}[Defi]{Theorem}

\newtheorem{exa}[Defi]{Example}
\newtheorem{Not}[Defi]{Notation}
\newcommand{\dR}{{\rm\scriptscriptstyle dR}}
\newcommand{\rk}{\mathop{\rm rank}\nolimits}
\newcommand{\Ind}{\mathop{\rm Ind}}
\makeatletter
\renewcommand{\subsection}{\@startsection{subsection}{2}%
        {\z@}{-3.25ex plus -1ex minus-.2ex}{-1em}{\bf}}
\makeatother

\begin{document}{\large}

\title{Teichm\"uller curves, triangle groups, and
 Lyapunov exponents}
\begin{abstract} 
  We construct a Teichm\"uller curve uniformized by the Fuchsian
  triangle group $\Delta(m,n,\infty)$ for every $m<n\leq \infty$.  Our
  construction includes the Teichm\"uller curves constructed by Veech
  and Ward as special cases. The construction essentially relies on
  properties of hypergeometric differential operators.  For small $m$,
  we find billiard tables that generate these Teichm\"uller
  curves.  We interprete some of the so-called Lyapunov exponents of
  the Kontsevich--Zorich cocycle as normalized degrees of a natural
  line bundle on a Teichm\"uller curve.  We determine the Lyapunov
  exponents for the Teichm\"uller curves we construct.
\end{abstract}
\date{\today}
\author{Irene I.~Bouw and Martin M\"oller}
\maketitle

\section*{Introduction} \label{Intro}
Let $C$ be a smooth curve defined over $\CC$. The curve $C$ is a {\em
  Teichm\"uller curve} if there exists a generically injective,
  holomorphic map from $C$ to the moduli space $M_g$ of curves of
  genus $g$ which is geodesic for the Teichm\"uller metric. Consider a
  pair $(X, \omega_X)$, where $X$ is a Riemann surface of genus $g$
  and $\omega_X$ is a holomorphic $1$-form on $X$. If the projective
  affine group, $\Gamma$, of $(X, \omega_X)$ is a lattice in
  $\PSL_2(\RR)$ then $C:=\HH/\Gamma$ is a Teichm\"uller curve. Such a
  pair $(X, \omega_X)$ is called a {\em Veech surface}. Moreover, the
  curve $X$ is a fiber of the family of curves $\cXX$ corresponding to
  the map $C\to M_g$. We refer to \S \ref{Teich} for precise
  definitions and more details.

Teichm\"uller curves arise naturally in the study of dynamics of billiard
paths on a polygon in $\RR^2$. Veech
(\cite{Ve89}) constructed a first class of Teichm\"uller curves $C=C_n$
starting from a triangle. The corresponding projective affine group
is commensurable to the triangle group $\Delta(2, n,  \infty)$. Ward
(\cite{Wa98}) also found  triangles which generate Teichm\"uller curves, with 
projective affine group   $\Delta(3,n,\infty)$. Several authors tried to find
other triangles which generate Teichm\"uller curves, but only sporadic
examples where found. Many types of triangles were disproven to
yield Veech surfaces (\cite{Vo96},(\cite{KeSm00}, \cite{Pu01}). 

In this paper we show that essentially all triangle groups
$\Delta(m,n,\infty)$ occur as the projective affine group of a
Teichm\"uller curve $C(m,n,\infty)$. (Since Teichm\"uller curves are
never complete (\cite{Ve89}), triangle groups $\Delta(m,n,k)$ with
$k\neq \infty$ do not occur.) We use a different construction from
previous authors; we construct the family $\cXX$ of curves defined by
$C$ rather than the individual Veech surface (which is a fiber of
$\cXX$). However, starting from our description, we compute an
algebraic equation for the corresponding Veech surface. The family
$\cXX$ is given as the quotient of an abelian cover ${\mathcal Y}\to
\PP^1$ by a finite group.

Under the simplifying assumption that $m<n<\infty$ and $n$ is odd, we
 relate the Veech surface corresponding to the Teichm\"uller curve
 $C(m,n,\infty)$ to a rational polygon. This polygon has $(m+3)/2$
 edges if $m$ if odd and $(m+4)/2$ edges if $m$ is even. This polygon
 does not have self-crossings if and only if $m\leq 5$. Therefore, for
 $m\leq 5$ we obtain  the Veech surface by unfolding a polygon.

From our construction we obtain new information even for the
Teichm\"uller curves found by Veech and Ward. Namely, we determine the
complete decomposition of the relative de Rham cohomology $R^1 f_* \CC_\cXX$
and the Lyapunov exponents, see below.

There exist Teichm\"uller curves whose projective affine group is not a
triangle group. McMullen (\cite{McM03}) constructed a series of such examples
in genus $g=2$. It would be interesting to try and extend our method to other
Fuchsian groups than triangle groups. This would probably be much more
involved due to the appearance of so-called {\em accessory parameters}.

We now give a more detailed description of our results.  Suppose that $m\geq
4$ and $m<n\leq \infty$ or that $m\geq 2$ and $3\leq n< \infty$.  We consider
a family of $N$-cyclic covers
$$
\cYY_t:\qquad  y^{N} = x^{a_1}(x-1)^{a_2}(x-t)^{a_3} $$
of the projective line
branched at $4$ points. Note that $\cYY$ defines a family over
$C=\PP^1_t-\{0,1,\infty\}$.  It is easy to  compute the differential equation
corresponding to the eigenspaces $\LL(i)$ of the action of $\ZZ/N$ on the
relative de Rham cohomology of $\cYY$ (\S \ref{covofPP}). These
eigenspaces are local systems of rank $2$, and the corresponding differential
equation is hypergeometric.  Cohen and Wolfart (\cite{CoWo90}) showed that
we may choose $N$ and $a_i$ in terms of $n$ and $m$ such that the projective
monodromy group of at least one of the eigenspaces $\LL(i)$ is the triangle
group $\Delta(m,n,\infty)$.

First consider the case that $m$ and $n$ are finite and relatively
prime.  Here we show that the particular choice of $N$ and the $a_i$
implies that, after replacing $C$ by a finite unramified cover, the
automorphism group of $\cYY$ contains a subgroup isomorphic to
$\ZZ/N\rtimes H$, where $H\simeq\ZZ/2\times \ZZ/2$. If $n$ is infinite
the group $H$ has order $2$. This case corresponds to half of Veech's
series of Teichm\"uller curves (\S \ref{ngonrev}).  If $m$ and
$n$ are not relatively prime we replace $\cYY$ by a suitable
$G_0$-Galois cover of the projective line, where $G_0$ is some
subgroup of $\ZZ/N\times \ZZ/N$. The description of $\cYY$ in this
case is just as explicit (\S \ref{realize}). 

 \par {\bf Theorem \ref{Teichneven} and \ref{hisTeich}:} {\em The
quotient family $\cXX:=\cYY/H$ is the pullback to $C$ of the universal family
over the moduli space of curves. The curve $C$ is an unramified
cover of a Teichm\"uller curve.}

\par The proof of this result relies on a Hodge-theoretical
characterization of Teichm\"uller curves (\cite{Mo04a}).  Another key
ingredient of the proof is the characterization of the vanishing of
the Kodaira--Spencer map in terms of invariants of the hypergeometric
differential equation corresponding to $\LL({i_0})$ (Proposition
\ref{KSvanunipotent}). Here $i_0$ is chosen such that the projective
monodromy group of $\LL({i_0})$ is the triangle group
$\Delta(m,n,\infty)$. The statement on the Kodaira--Spencer map
translates to the following geometric property of $\cXX$. A fiber
$\cXX_c$ of $\cXX$ is singular if and only if the monodromy around $c$
of the local system induced by $\LL({i_0})$ is infinite (Proposition
\ref{Higgsprop}). This is one of the central observations of the
paper. This is already apparent in our treatment of the relatively
straigthforward case of Veech's families of Teichm\"uller curves in
\S \ref{ngonrev}.


{\bf Theorem \ref{mnwhichtriangle}:} {\em Suppose that $n$ is finite and 
  $m$ is different from $n$. Then the projective affine group of $\cXX$ is the
  triangle group $\Delta(m,n, \infty)$.}

We determine the projective affine group of our Teichm\"uller curves
directly from the construction of the family $\cXX$ and do not need to
consider the corresponding Veech surfaces, as is done by Veech and
Ward. For example, we  determine the number of zeros of the
generating differential of a Veech surface corresponding to
$C(m,n,\infty)$ in terms of $n$ and $m$ by algebraic methods (Theorem
\ref{signmninf}).


In \S \ref{Billiards} we change perspective, and discuss the
question of realizing our Teichm\"uller curves via unfolding of
rational polygons (or: billiard tables). This section may be read
independently of the rest of the paper. For $m\leq 5$ we construct a
billiard table $T(m,n,\infty)$ and show that it defines a
Teichm\"uller curve, via unfolding. For $m=2,3$ this gives the
triangles considered by Veech (\cite{Ve89}) and Ward
(\cite{Wa98}). For $m=4,5$ we find new billiard tables which are
rational $4$-gons. We interpret the Veech surfaces corresponding to
these billiard tables as  fiber 
 of the family $\cXX \to C$ of curves.
A key ingredient here is a theorem of Ward (\cite{Wa98},
Theorem C') which relates a cyclic cover of the projective line to a
polygon, via the Schwarz--Christoffel map. We then use that certain
fibers of $\cXX$ are a cyclic cover of the projective line (Theorem
\ref{Wardexplicit}).

For $m\geq 6$  the same procedure still produces rational 
polygons $T(m,n,\infty)$, but they have self-crossings
and therefore do not define billiard tables. In principle, one 
could still describe the translation surface corresponding to
$T(m,n,\infty)$, but these would be hard to visualize.


Our last main result concerns Lyapunov exponents. Let $V$ be a flat
normed vector bundle on a manifold with flow. The {\em Lyapunov
exponents} measure the rate of growth of the length of vectors in $V$
under parallel transport along the flow. We refer to \S
\ref{Lyapunov} for precise definitions and a motivation of the
concept. We express the Lyapunov exponents for an arbitrary
Teichm\"uller curves in terms of the degree of certain local systems.
\par Let $f:\cXX\to C$ be the universal family over an unramified
cover of an arbitrary Teichm\"uller curve. The relative de Rham
cohomology $R^1 f_\ast \CC_\cXX$ has  $r$ local subsystems $\LL(i)$
of rank two. The associated vector bundles carry a Hodge filtration
(Theorem \ref{maxHiggsTeich}).  The $(1,0)$-parts of the Hodge
filtration are line bundles $\cLL(i)$ and the ratios
$$\lambda_i := 2\deg(\cLL(i)) /(2g(C)-2+s), \quad s = {\rm card}(\ol{C} \sms C)
$$
are unchanged if we pass to an unramified cover of $C$.
\par
{\bf Theorem \ref{rLyap}:} {\em The ratios $\lambda_i$ are 
$r$ of $g$ non-negative Lyapunov exponents of the 
Kont\-sevich--Zorich cocycle over the Teich\-m\"uller geodesic flow
on the canonical lift of 
a Teich\-m\"uller curve to the $1$-form bundle over the
moduli space.}
\par

A sketch of the relation between the degree of
$f_* \omega_{\cXX/C}$ and the sum of all Lyapunov exponents
already appears in \cite{Ko97}.

Now suppose that $C$ is an unramified cover of $C(m,n,\infty)$ (
Theorems \ref{Teichneven} and \ref{hisTeich}), and let $f:\cXX\to C$
be the corresponding family of curves.  In Corollaries
\ref{LyapnevenXX} (Veech's series), \ref{isVeechneven} and
\ref{splitVHSh} we give an explicit expression for all Lyapunov
exponents of $C$. For Veech's series of Teichm\"uller curves and for a
series of square-tiled coverings the Lyapunov exponents were
calculated independently by Kontsevich and Zorich (unpublished). They
form an arithmetic progression in these cases.  Example \ref{lyapex}
shows that this does not hold in general.  
\par 
It is well-known that
the largest Lyapunov exponent $\lambda_1=1$ occurs with multiplicity
one.  We interpret $1-\lambda_i$ as the number of zeros of the
Kodaira--Spencer map of $\LL(i)$, counted with multiplicity (\S
\ref{Teich}), up to a factor. For the Teichm\"uller curves constructed
in Theorems \ref{Teichneven} and \ref{hisTeich} we determine the
position of the zeros of the Kodaira--Spencer map.  These zeros are
related to elliptic fixed points of the projective affine group
$\Gamma$ (Propositions \ref{KSvanunipotent} and
\ref{lambdalocsyst}). For an arbitrary Teichm\"uller curve it is an
interesting question to determine the position of the zeros of the
Kodaira--Spencer map. Precise information on the zeros of the
Kodaira--Spencer map might shed new light on the defects $1-\lambda_i$ of the
Lyapunov exponents.  \par The starting point of this
paper was a discussion with Pascal Hubert and Anton Zorich on Lyapunov
exponents.  The second named author thanks them heartily. Both authors
acknowledge support from the DFG-Schwerpunkt `Komplexe
Mannigfaltigkeiten'.  We thank Frits Beukers for suggesting the proof
of Proposition \ref{Betaprop} and  Silke Notheis for
cartographic support.  \par

\section{Teichm\"uller curves} \label{Teich}

A {\em Teichm\"uller curve} is a generically injective, 
holomorphic map $C \to M_g$ from a smooth algebraic curve 
$C$ to the moduli space of curves of genus $g$ which is
geodesic for the Teichm\"uller metric. A Teichm\"uller curve 
arises as quotient $C = \HH/\Gamma$,
where $\HH \to T_g$ is a complex Teichm\"uller geodesic in Teichm\"uller
space $T_g$. Here $\Gamma$ is the subgroup in the Teichm\"uller
modular group fixing $\HH$ as a subset  of $T_g$ (setwise, not pointwise)
and  $C$ is the normalization of the image $\HH \to T_g \to M_g$.

Veech showed that a Teichm\"uller curve $C$ is never complete
(\cite{Ve89} Prop.\ 2.4). We let $\ol{C}$ be a smooth completion of
$C$ and $S := \ol{C} \sms C$.  In the sequel, rather than consider
Teichm\"uller curves themselver,  it will be convenient to
consider  finite unramified
covers of $C$ that satisfy two conditions: the corresponding
subgroup of $\Gamma$ is torsion free and the moduli map factors
through a fine moduli space of curves (e.g.\ with level structure
$M_g^{[n]}$).  We nevertheless stick to the notation $C$ for the base
curve and let $f:\cXX \to C$ be the pullback of the universal family
over $M_g^{[n]}$ to $C$. We will use $\ol{f}: \ol{\cXX} \to \ol{C}$
for the family of stable curves extending $f$. See also \cite{Mo04a}
\S 1.3.  \par Teichm\"uller curves, or more generally geodesic
discs in Teichm\"uller space, are generated by a pair $(X, q)$ of a
Riemann surface and a quadratic differential $q \in
\Gamma(X,(\Omega^1_X)^{\otimes 2})$.  These pairs are called {\em
translation surfaces}.  If a pair $(X,q)$ generates a Teichm\"uller
curve, the pair is called a {\em Veech surface}. Any smooth fiber of
$f$ together with the suitable quadratic differential is a Veech
surface.  Theorem \ref{maxHiggsTeich} below characterizes
Teichm\"uller curves where $q = \omega^2$ is the square of a
holomorphic $1$-form $\omega \in \Gamma(X,(\Omega^1_X))$. The examples
we construct will have this property, too. Hence: \par {\em From now
on the notion `Teichm\"uller curve' includes `generated by a
$1$-form'.}  \par For a pair $(X,\omega)$ we let ${\mathop{\rm
Aff}}^+(X,\omega)$ be the group of orientation preserving
diffeomorphism of $X$ that are affine with respect to the charts
provided by integrating $\omega$.  Associating to an element of
${\mathop{\rm Aff}}^+(X,\omega)$ its matrix part gives a well-defined
map to $\SL_2(\RR)$.  The image of this map in $\SL(X,\omega)$ is
called the {\em affine group} of $(X,\omega)$.  The matrix part of an
element of ${\mathop{\rm Aff}}^+(X,\omega)$ is also called its {\em
derivative}.  The stabilizer group $\Gamma$ of $\HH \hookrightarrow
T_g$ coincides, up to conjugation, with the affine group
$\SL(X,\omega)$ (\cite{McM03}). We denote throughout by $K =
\QQ({\rm tr}(\gamma, \gamma \in \Gamma))$ the trace field and let $r
:= [K:\QQ]$. We call the image of $\SL(X,\omega)$ in $\PSL_2(\RR)$ the
{\em projective affine group} and denote it by $\PSL(X,\omega)$.
\newline
We refer to \cite{KMS86} and \cite{KeSm00} for a systematic
description of Teichm\"uller curves in terms of billiards.  
\par 
We
recall from \cite{Mo04a} Theorem~2.6 and Theorem~ 5.5 a description of
the variation of Hodge structures (VHS) over a Teichm\"uller curve,
and a characterization of Teichm\"uller curves in these terms.
\newline
Let $\LL$ be a rank two irreducible $\CC$-local system on an affine curve $C$.
Suppose that the Deligne extension $\cEE$ of $\LL \otimes_\CC \cOO$ 
(\cite{De70} Proposition II.5.2) to $\ol{C}$ carries
a Hodge filtration of weight one $\cLL:= \cEE^{(1,0)} \subset \cLL$. 
We denote by $\nabla$ the corresponding logarithmic connection on $\cEE$.
The {\em Kodaira--Spencer map} (also:
Higgs field, or:  second fundamental form) with respect to $S$ is
the composition map
\begin{equation} \label{KSmap}
\Theta: \cLL \to \cEE \stackrel{\nabla}{\to}
 \cEE \otimes \Omega^1_{\ol{C}}(\log S) 
\to (\cEE/\cLL) \otimes \Omega^1_{\ol{C}}(\log S). 
\end{equation}
A  VHS of rank $2$ and weight one whose Kodaira--Spencer map
with respect to some $S$ vanishes nowhere on $\ol{C}$ is called {\em
maximal Higgs} in \cite{ViZu04}. The corresponding vector
bundle $\cEE$ is called {\em indigenous bundle}. See \cite{BoWe05}
or \cite{Mo99} for appearances of such bundles with more emphasis
on char $p > 0$. 
\par
\begin{Thm} \label{maxHiggsTeich}
\begin{itemize}
\item[(a)]
Let $f: \cXX \to C$ be the universal family over a finite unramified
cover of a Teichm\"uller curve. Then we have a decomposition of
the VHS of $f$ as
\begin{equation} \label{VHSdecomp}
 R^1 f_* \QQ = \WW \oplus \MM \quad \text{and} \quad
\WW \otimes_ \QQ \CC  = \bigoplus_{i=1}^r \LL_i.
\end{equation}
In this decomposition the $\LL_i$ are Galois conjugate, irreducible,
pairwise non-isomorphic, $\CC$-local systems of rank two.
The $\LL_i$ are in fact defined over some field $F \subset \RR$ 
that is Galois over $\QQ$ and contains the trace field $K$. 
Moreover, $\LL_1$ is maximal Higgs.
\item[(b)]
Conversely, suppose $f: \cXX \to C$ is a family of smooth curves
such that $R^1 f_* \CC$ contains a local system of rank two 
which is maximal Higgs with respect to the set $S = \ol{C} \sms C$. 
Then $f$ is the universal family over a finite unramified
cover of a Teichm\"uller curve.
\end{itemize}
\end{Thm}
\par
Note that `maximal Higgs'  depends on $S$. We will encounter cases 
where $\LL$ extends over some points of $S$ and becomes maximal
Higgs with respect to a smaller set $\Su \subset S$, but it
is not maximal Higgs with respect to  $S$. See also
Proposition \ref{Higgsprop} and Remark \ref{allparrem}.

\section{Local exponents of differential equations
and zeros of the Kodaira--Spencer map} \label{locexpo}

In this section we provide a dictionary between local systems plus a
section on the one side and differential equations on the other
side. In particular, we translate local properties of a differential
operator into vanishing  of the Kodaira--Spencer map. In the \S\S 
\ref{ngonrev} and \ref{realize} we essentially start
with a hypergeometric differential equation whose local properties are
well-known. Via Proposition \ref{KSvanunipotent} the vanishing of the
Kodaira--Spencer map of the corresponding local system is completely
determined.  This knowledge is then exploited in a criterion
(Proposition \ref{Higgsprop}) for a family of curves $f:\cXX \to C$ to
be the universal family over a Teichm\"uller curve.  \par Let $\LL$ be
a irreducible $\CC$-local system of rank $2$ on an affine curve $C$,
not necessarily a Teichm\"uller curve. Let $C\hookrightarrow \ol{C}$ be
the corresponding complete curve, and let $\cEE$ be the Deligne
extension of $\LL$ (\S \ref{Teich}).  We suppose that $\LL$
carries a polarized VHS of weight one and choose a section $s$ of
$(\LL \otimes_\CC \cOO_C)^{(1,0)}$.  Let $t$ be a coordinate on $C$.
We denote by $D:= \nabla({\partial}/{\partial t})$. Since $\LL$ is
irreducible, the sections $s$ and $Ds$ are linearly independent.
Hence $s$ satisfies a differential equation $Ls=0$, where
$$L = D^2 + p(t) D + q(t),$$ for some meromorphic functions $p,q$ on
$\ol{C}$. Note that we may interprete $L$ as a second order
differential operator $L:\cOO_C\to\cOO_C$, by interpreting $D$ as
derivation with respect to $t$.  \par Conversely, the set of solutions
of a second order differential operator $L: \cOO_C \to \cOO_C$ forms a
local system ${\rm Sol} \subset\cOO_C$. If $L$ is obtained from $\LL$
then ${\rm Sol} \cong \LL^\vee$ (\cite{De70} \S 1.4).  The canonical
map
$$\varphi: {\rm Sol}\, \otimes_\CC \cOO_C \to \cOO_C, \quad f \otimes g
\mapsto fg$$
hence defines a section $s =s_\varphi$ of $\LL \otimes_\CC
\cOO_C$.  \par A point $c \in \ol{C}$ is a {\em singular point} of $L$ if $p$
or $q$ has a pole at $c$.  In what follows, we always assume that $L$ has
regular singularities.  Let $t$ be a local parameter at $c\in \ol{C}$. Recall
that $L$ has a {\em regular singularity} at $c$ if $(t-c)p$ and $(t-c)^2q$ are
holomorphic at $c$, by Fuchs' Theorem. Note that there is a difference between
the notions `singularity of the Deligne extension of the local system $\LL$'
and `singularities of the differential operator $L$'. We refer to \cite{Ka70},
\S 11 for a definition of the notion regular singularity of a flat vector
bundle. (The essential difference between the two notions is that the basis
${\mathbb e}$ of \cite{Ka70} (11.2.1), need not be a cyclic basis
(\cite{Ka70} \S 11.4).) 
Unless stated explicitly, we only use the notion of singularity of the
differential operator.

  The {\em local
exponents} $\gamma_0$, $\gamma_1$ of $L$ at $c$ are the roots of the
characteristic equation
$$ t(t-1) + tp_{-1} +q_{-2} =0, $$
where $p = \sum_{i=-1}^\infty p_i (t-c)^i$ and 
$q = \sum_{i=-2}^\infty q_i (t-c)^i$. The table
recording singularities and the local exponents is usually 
called {\em Riemann scheme}. See e.g.\ \cite{Yo87} \S 2.5 for more details.
\newline
Note that $L$ and the local exponents not only depend on $\LL$ but
also on the section chosen. Replacing $s$ by $\alpha s$ shifts the
local exponents at $c$ by the order of the function $\alpha$ at
$c$. The exponentials $e^{2\pi it_1}$ and $e^{2\pi it_2}$ of the
local exponents are the eigenvalues of the local monodromy
matrix of $L$ at $c$.
 The following criterion
is well-known (e.g.\ \cite{Yo87} \S I.2.6).
\par
\begin{Lemma} \label{unipotent}
All local monodromy matrices  of ${\rm Sol}$ are unipotent
if and only if  both local exponents are
integers for all $c\in \ol{C}$.
\end{Lemma}
\par In the classical case that $\ol{C} \cong \PP^1$ the differential
operator $L$ is determined by the local exponents exactly if the
number of singularities is three; this is the case of hypergeometric
differential equations. We will exploit this fact in the next
sections. If the number of singularities is larger than three, $L$ is
no longer determined by the local exponents and the position of the
singularities, but also depends on the {\em accessory parameters}
(\cite{Yo87} \S I.3.2).

In the rest of this section we suppose that all local monodromy
matrices of $L$ are unipotent. We define $S_{\scriptstyle\rm
u}=\ol{C}-C$ as the set of points where the monodromy is
nontrivial. Let $S \subset \overline{C}$ be a set containing the
singularities of the Deligne extension of $\LL$. The reader should
think of $S$ being the set of singular fibers of a family of
curves over $\ol{C}$. In particular $S \supset \Su$.

The following proposition expresses the order of vanishing of the
Kodaira--Spencer map $(\ref{KSmap})$ at $c \in \ol{C}$ in terms of 
the local exponents at $c$. If $c \in C$ we suppose that the section $s$ 
is chosen such that the local exponents are $(0, n_c)$ with $n_c\geq 0$.
This is always possible, multiplying $s$ with a power of a local
parameter if necessary.  
\par
\begin{Prop} \label{KSvanunipotent}
\begin{itemize}
\item[(a)] Let $c\in C$. Then  $n_c\geq 1$.
\item[(b)] Suppose that $c \not\in S$.  The order of vanishing of $\Theta$ at 
$b$ is  $n_c-1$.
\item[(c)] Suppose that $c \in S\sms S_{\scriptstyle\rm u}$.  The order of
 vanishing of $\Theta$ at $b$ is $n_c$.
\item[(d)] If $c \in \Su$
then $\Theta$ does not vanish at $c$. 
\end{itemize}
\end{Prop}
\par {\bf Proof:} Suppose that $c \in C$. Our assumptions imply that
the local exponents $(0, n_c)$ at $c$ are nonnegative integers. Since
$\LL$ is a local system on $C$, it has two linearly independent
algebraic section in a neighborhood of $c$. This implies that $n_c\geq
1$ (\cite{Yo87} \S I.2.5). This proves (a).

 If $c \not\in S$ the differential operator $L$ has 
solutions $s_1, s_2$ with leading terms $1$ and $t^{n_c}$,
respectively (\cite{Yo87} I, 2.5).  We want to determine the vanishing
order of $D(s)$ in $\cEE/(s \otimes_\CC \cOO_C)$. By the above
correspondence between the local system and the differential equation
we may as well calculate the vanishing order of $D(\varphi)$ in $({\rm
Sol}^\vee \otimes_\CC \cOO_C)/ (\varphi \otimes_\CC \cOO_C)$.
A basis of ${\rm Sol}^\vee \otimes_\CC \cOO_C$ around $c$ is 
$$s_i^\vee:\quad s_1 \otimes g_1 + s_2 \otimes g_2 \mapsto
s_ig_i \quad (i=1,2).$$
By definition of the dual connection and the flatness of $s_i$
one calculates that $D(\varphi)$ is the class of 
$$ s_1 \otimes g_1 + s_2 \otimes g_2 \mapsto g_1 s'_1 + g_2 s'_2$$ in
$({\rm Sol}^\vee \otimes_\CC \cOO_C)/ (\varphi \otimes_\CC \cOO_C)$.
Since both $\varphi$ and $s_1$ do not vanish at $c$, we conclude
that the order of vanishing of $D(\varphi)$ at $c$ is  $n-1$. This proves (b).
\newline
In the  case that $c \in S$ we should consider the contraction against $t
{\partial}/\partial t$. This increases the order of vanishing
of $\Theta$ by one. This proves (c).  \par We now treat the case that
$c \in \ol{C} \sms C$. Consider the residue map ${\rm Res}_c(\nabla)
\in {\rm End}(\cEE_c)$. Suppose the Kodaira--Spencer map vanishes at
$c$. This implies that ${\rm Res}_c(\nabla)$ is a diagonal matrix in a
basis consisting of an element from $\cLL_c$ and an element from its
orthogonal complement. But ${\rm Res}_c(\nabla)$ is nilpotent
(\cite{De87} Proposition II.5.4 (iv)), hence zero.  This implies that
two linearly independent sections of $\LL$ extend to $c$. This
contradicts the hypothesis on the monodromy around $c$. This proves
(d).  
\hspace*{\fill} $\Box$ 
\par 
The ratios $\lambda(\LL, S) :=
2\deg(\cLL)/\Omega^1_{\ol{C}}(\log S)$ will be of central interest in
the sequel. The factor $2$ is motivated by  \S \ref{Lyapunov},
where we interprete the $\lambda(\LL, S)$ as Lyapunov
exponents. Therefore we call the $\lambda(\LL, S)$ from now on {\em
Lyapunov exponents}. We will suppress $S$ if it is clear from the
context.  \par
\begin{Rem} \label{VHSunique} {\rm
We will only be interested in local $\CC$-systems $\LL$ that arise as
local subsystems of $R^1 f_* \CC$ for a family of curves $f: \cXX \to
C$. In this case a Hodge filtration exists on $\LL$ and is unique
(\cite{De87} Prop.\ 1.13). Therefore we only have to keep track of the
local system, but not of  the VHS.
}\end{Rem} 
\par
 The following
lemma is noted for future reference. The proof of straightforward.  \par
\begin{Lemma} \label{quotindepent}
The ratio $\lambda(\LL, S)$
does not change by taking unramified coverings.
\end{Lemma}
\par

\section{Cyclic covers of 
the projective line branched at $4$ points} \label{covofPP}

Let $N>1$ be an integer, and suppose given a $4$-tuple of integers $(a_1,
\ldots, a_4)$ with $0<a_\mu<N$ and $\sum_{\mu=1}^4 a_\mu=(k+1)N$, for some
integer $k$.  We denote by $\PP^1$ the projective line with coordinate $t$,
and put  $\PP^*=\PP^1-\{0,1,\infty\}$. Let $\cPP\simeq\PP^1 \times \PP^* 
\to \PP^*$ be the trivial fibration with fiber coordinate $x$.  
Let $x_1=0, x_2=1, x_3=t, x_4=\infty$ be
sections of $\cPP \to \PP^*$.  We fix an injective character $\chi: \ZZ/N
\to \CC^*$.  Let $g: \cZZ \to \PP^*$ be the $N$-cyclic cover of type $(x_\mu,
a_\mu)$ (\cite{Bo04} Definition 2.1). This means that $\cZZ$ is the family of
projective curves with affine model
\begin{equation} \label{nfoldcov}
\cZZ_t: \qquad z^N = x^{a_1}(x-1)^{a_2}(x- t)^{a_3}.
\end{equation}
 We suppose, furthermore, that
$\gcd(a_1, a_2, a_3, a_4, N)=1$. This implies that the family is connected.
The genus of $\cZZ_t$ is $N+1-(\sum_{\mu=1}^4\gcd(a_\mu, N))/2$.
\par
In this section, we collect some well-known  facts on such cyclic covers.
We write 
\[
\sigma_\mu(i) = \langle i a_\mu/N \rangle=a_\mu(i)/N,
\]
where $\langle \cdot \rangle$ denotes the fractional part.  Let $k(i)+1 =
\sum_{\mu=1}^4 \sigma_\mu(i) $. We fix an injective character $\chi: \ZZ/N \to
\CC^*$ such that $h \in \Gal(\cZZ/\cPP) \cong \ZZ/N$ acts as $h \cdot z=
\chi(h)z$.
\par
\begin{Lemma}\label{dimlem}
 For $0<i<N$, we let\,  $s(i)$ be the number of\,  $a_\mu$ unequal to\,  $0
 \bmod{N/\gcd(i, N)}$. Put  $\LL(i)= H^1_\dR (\cZZ/\PP^*)$. Then
\begin{itemize}
\item[(a)] $\dim_\CC\LL(i)=s(i)-2$,
\item[(b)] $\rk g_* (\Omega^1_{\cZZ/\PP^*})_{\chi^i}=s(i)- 2-k(i), \qquad
  \rk (R^1g_* \cOO_\cZZ)_{\chi^i}=k(i).$
\item[(c)] If $k(i)=1$ then
\[
 \omega_i := \frac{z^i\,{\rm d}x}{x^{1+[i\sigma_1]}
(x-1)^{1+[i\sigma_2]}(x-t)^{1+[i\sigma_3]}}
\]
is a non-vanishing section of
$g_* (\Omega^1_{\cZZ/\PP^*})_{\chi^i}$. It  is a solution 
of the hypergeometric differential operator
\[
L(i) := \nabla\!\left(\frac{\partial}{\partial t} \right)^2 +
\frac{(A(i)+B(i)+1)t -C(i)}{t(t-1)} \,\,\nabla\!\left(\frac{\partial}
{\partial t} \right)
+ \frac{A(i)B(i)}{t(t-1)}, 
\]
where $ A(i) = 1-\sigma_3(i), \quad B(i)=
2-(\sigma_1(i)+\sigma_2(i)+\sigma_3(i)), \quad C(i) = 2-(\sigma_1(i) +
\sigma_3(i)).$
\end{itemize}
\end{Lemma}
{\bf Proof:} The second statement of (b) is proved in \cite{Bo01}
Lemma 4.3. The first statement follows from Serre duality and
\cite{Bo01} Lemma 4.5. Part (a) follows immediately from (b).  The
statement that $\omega_i$ is holomorphic and non-vanishing is a
straightforward verification. The statement that $L(i) \omega_i=0$ in
$H^1_\dR(\cZZ/\PP^*)_{\chi^i}$ is proved for example in \cite{Bo05},
Lemma 1.1.4. \hspace*{\fill} $\Box$ \par The differential operator $L(i)$
corresponds to the local system $\LL(i)=H^1_\dR(\cZZ/\PP^*)_{\chi^i}$
together with the choice of a section $\omega_i$ via the
correspondence described at the beginning of \S \ref{locexpo}. It
has singularities precisely at $0$, $1$ and $\infty$. 
Its local exponents are
summarized in the Riemann scheme
\begin{equation} \label{RSHG}
\left[ \begin{array}{ccc}
t=0&t=1&t=\infty \\
0  & 0 & A(i) \\
\gamma_0:=1-C(i) &  \gamma_1:=C(i) -A(i)-B(i) & B(i) \\
\end{array}\right].
\end{equation}
\par
A (Fuchsian) $(m,n,p)$-{\em triangle group} for $m,n,p \in \NN \cup \{\infty\}$
satisfying $1/m+1/n+1/p < 1$ is a Fuchsian group in $\PSL_2(\RR)$
generated by matrices $M_1, M_2, M_3$ satisfying $M_1M_2M_3 = 1$ and
$$ \tr(M_1) = \pm 2\cos(\pi/m), \quad \tr(M_2) = \pm 2\cos(\pi/n),
\quad \tr(M_3) = \pm 2\cos(\pi/p). \quad$$ A triangle group is
determined, up to conjugation in $\PSL_2(\RR)$, by the triple
$(m,n,p)$.  It is well-known that the projective monodromy groups of
the hypergeometric differential operators $L(i)$ are triangle groups
under suitable conditions on $A(i), B(i), C(i)$. These conditions are
met in the cases we consider in \S \ref{ngonrev} and
\ref{realize}.  \par We are interested in determining the order of
vanishing of the Kodaira--Spencer map. Note that if $k(i)=0$ or
$k(i)=2$ then the Hodge filtration on the corresponding eigenspace is
trivial and hence the Kodaira--Spencer map is zero.  \par Let
$\overline{\pi}:\overline{C}\to \PP^1$ a finite cover, unbranched 
outside $\{0,
1, \infty\}$, such that the monodromy of the pullback of $\cZZ$ via
$\overline{\pi}$ is unipotent for all $c\in \overline{C}$.  

Let $\Su=\Su(i)\subset \pi^{-1}(0,1,\infty)$ be the set of points such
that $\LL(i)$ has nontrivial local monodromy. Our assumption implies
that the monodromy at $c\in \Su$ is infinite. In what follows, the set
$\Su$ will be nonempty. It is therefore no restriction to suppose that
$\pi^{-1}(\infty)$ is contained in $\Su$.  In terms of the invariants
$a_\mu$ this means that $a_3(i)+a_4(i)\equiv 0\bmod{N}$. It follows
that $A(i)=B(i)$.  Let $b_0$ (resp.\ $b_1$) be the common denominator of
the local exponents $\gamma_0(i)$ (resp.\ $\gamma_1(i)$) for $1 \leq i
< N$. Write $|\gamma_0(i)| = n_0(i)/b_0$ and
$|\gamma_1(i)|=n_1(i)/b_1$.  Note that $\pi^{-1}(t=\mu)\subset
\Su(i)$ if and only if $\gamma_\mu(i)=0$. Therefore
$\Su(i)=\pi^{-1}(\{0,1,\infty\})$ if and only if
$\gamma_0(i)=\gamma_1(i)=0$ . It easily follows that the set $\Su$ is
in fact independent of $i$.


The following proposition is
the basic criterion we use for constructing Teichm\"uller curves. 
\par
\begin{Prop} \label{Higgsprop} 
Consider a family of curves $\cZZ\to \PP_t^\ast$ as in
(\ref{nfoldcov}). Let $ 0<i_0<N$ be an integer such that
\begin{equation}\label{Higgseq}
\gamma_\mu(i_0)=1/b_\mu\quad\mbox{for all }\mu\in \{0,1\}\mbox{ with }
\gamma_\mu(i_0)\neq 0.
\end{equation}
There is a finite cover $\pi:\ol{C}\to \PP^1_t$ branched of order exactly
$b_\mu$ at $t=\mu \in \{0,1\}$ for all $\mu$ such that $\gamma_\mu(i_0)\neq
0$. Moreover, we require that the local monodromy of the pullback of
$\LL(i_0)$ to $\ol{C}$ is unipotent, for all $c\in \ol{C}$. Write $\cZZ_{{C}}$
for the pullback of $\cZZ$ to ${C}$.
\par
Choose a subgroup $H$ of $\Aut(\cZZ_C)$ 
and define $\cXX:=\cZZ_C/H$.  Suppose that
\begin{itemize}
\item $\cXX$ extends to a smooth family over 
$\widetilde{C} :=\ol{C} \sms \Su$,
\item  there is a local system $\LL$ isomorphic to $\LL(i_0)$ which
descends to $\cXX$.
\end{itemize} 
Then the moduli map $\widetilde{C} \to M_g$ is an unramified cover of
a Teichm\"uller curve.
\end{Prop}
\par
This criterion will be applied to subgroups $H$ that intersect
trivially with $\Gal(\cZZ_C/\cPP_C)$.
\par
{\bf Proof:} If $\gamma_\mu(i_0)\neq 0$ the monodromy of $g$ at
$t=\mu$ becomes trivial after pullback by a cover which is branched at
$t=\mu$ of order $b$ if and only if $b_\mu|b$. Hence if the cover
$\ol{\pi}$ is sufficiently branched at points over $\Su$, the local monodromy
of the pullback of $\LL(i_0)$ to $C$ is  unipotent by Lemma \ref{unipotent}.
\par 
The local exponents of the pullback of $\LL(i_0)$ to $\ol{C}$
are the original ones multiplied by the ramification index. Hence, for all
$c\in \widetilde{C}$,  the local exponents  
are $(0,1)$. By definition, the same holds for the local exponents of the
bundle $\LL$.
The hypothesis
on the singular fibers of $\cXX$ implies that $c\in \widetilde{C}$ is not a
singularity of the flat bundle $\LL$ (\cite{Ka70} \S 14).
Therefore we may apply Proposition \ref{KSvanunipotent} (b) and (d) to $\LL$
with $S=\Su$. We conclude that 
the Kodaira--Spencer map of $\LL$  vanishes nowhere.
The proposition therefore follows from Theorem 1.1. 
\hspace*{\fill} $\Box$
\par
\begin{Rem} \label{allparrem}{\rm 
The structure of the stable model $g_{\ol{C}}$ of the family $g_C:
\cZZ_C \to C$ is given in the next subsection.  It implies that all
fibers of preimages of $\{0,1,\infty\}$ are singular. Hence applying
Proposition \ref{Higgsprop} to $g_{\ol{C}}$ with $H=\{1\}$, we find
that $g_{\ol{C}}$ defines a Teichm\"uller curve if and only if
$\Su=\pi^{-1}(\{0,1,\infty\})$. This happens for example for the families
$$ y^2=x(x-1)(x-t) \quad \text{and} \quad y^4=x(x-1)(x-t).$$
Here $\overline{C}=\PP^1_t$, and the uniformizing group is the triangle group
$\Delta(\infty, \infty, \infty)$. Clearly, this is a very special
situation. 
}\end{Rem}
 
\begin{Prop} \label{lambdalocsyst}
  Let $0<i<N$ be an integer with $k(i)=1$.  Denote by $\cLL(i)$ the
  $(1,0)$-part of the local system $\LL(i)$ over $C$. 
Then
\[
 \deg \cLL(i) = \frac{\deg(\pi)}{2}\left(1-\frac{n_0(i)}{b_0}-
\frac{n_1(i)}{b_1}\right)
\]
with the convention that $1/b_\mu = 0$ if $n_\mu=0$. In particular, 
the Lyapunov exponent
$$ \lambda(\LL(i),\Su) = \left(1-\frac{n_0(i)}{b_0}-
\frac{n_1(i)}{b_1}\right) / \left(1-\frac{1}{b_0}-
\frac{1}{b_1}\right) $$
is independent of the choice of $\ol{\pi}$.
\end{Prop}
\par
{\bf Proof:}
We only treat the case that both $n_0(i)$ and
$n_1(i)$ are non-zero, leaving the few modifications in the other cases to
the reader. One checks that
$$ \deg \Omega^1_{\ol{C}}(\log \Su) = \deg(\ol{\pi})\left(1-\frac{1}{b_0}-
\frac{1}{b_1}\right)$$ is independent of the ramification order of $g$ over
$t=\infty$.  It follows from the definition (\ref{KSmap}) of the
Kodaira--Spencer map $\Theta$ that $2 \deg \cLL_{\chi^i}-\deg
\Omega^1_c(\log \Su)$ is the number of zeros of $\Theta$, counted with
multiplicity.  Therefore the proposition follows from Proposition
\ref{KSvanunipotent}. 
\hspace*{\fill} $\Box$ 

\subsection{Degenerations of cyclic covers} 
We now describe the stable model of the degenerate fibers of $\cZZ$. For
simplicity, we only describe the fiber $\cZZ_0$ above $t=0$. The other
degenerate fibers may be described similarly, by permuting $\{0,1,t,
\infty\}$. A general reference for this is \cite{Stefandiss} \S 4.3.
However, since we consider the easy situation of cyclic covers of the
projective line branched at $4$ points, we may simplify the presentation.
\par
As before, we let $\cPP\to \PP^*$ be the trivial fibration with fiber
coordinate $x$.
We consider the sections $x_1=0, x_2=1, x_3=t, x_4=\infty$ of $\cPP\to \PP^*$
as marking on $\cPP$. We may extend $\cPP$ to a family of stably marked curves
over $\PP(=\PP^1_t)$, which we still denote by $\cPP$. The fiber $P_0$ of
 $\cPP$ at
$t=0$ consists of two irreducible components which we denote by $P_0^1$ and
$P_0^2$. We assume that $x_1$ and $x_3$ (resp.\ $x_2$ and $x_4$) specialize
to the smooth part of $P_0^1$ (resp.\ $P_0^2$). We denote the
intersection point of $P_0^1$ and $P_0^2$ by $\xi$. It is well known
that the family of curves $f:\cZZ\to \cPP$ over $\PP^\ast$ extends to a family
of {\sl admissible covers} over $\PP^1_t$. See for example \cite{HaSt} or
\cite{We99}.  For a short overview we refer to \cite{An} \S 2.1.
\par
The definition of type (\cite{Bo04} Definition 2.1) immediately implies that
the restriction of the admissible cover $f_0:\cZZ_0\to \cPP_0$ to $P_0^1$
(resp.\ $P_0^2$) has type $(x_1, x_3, \xi; a_1, a_3, a_2+a_4)$ (resp.\ 
 $(x_2, x_4, \xi; a_2, a_4, a_1+a_3)$). (Admissibility amounts 
in our situation
to $(a_1+a_3)+(a_2+a_4)\equiv 0\bmod{N}$.) Let $Z_0^j$ be a connected
component of the restriction of $\cZZ_0$ to $P_0^j$. Choosing suitable
coordinates, $Z_0^1$ (resp. $Z_0^2$) is a connected component of the smooth
projective curve defined by the equation $z^N=x^{a_1}(x-1)^{a_3}$ (resp.\ 
 $z^N=x^{a_2}(x-1)^{a_4}$).

Denote by $H^j=\Gal(Z^j_0, P^j_0)\subset H\simeq \ZZ/N$ the subgroups
obtained by restricting the Galois action. 
Then $\cZZ_0$ is obtained by suitably identifying the
points in the fiber above $\xi$ of 
$\Ind_{H^1}^H Z_0^1 $ and $\Ind_{H^2}^H Z_0^2$. 
\par
Proposition \ref{Degeneration} follows from the explicit description of
the components of $\cZZ_0$.   Put $\beta_1=\gcd(a_1,
a_3, N)$ and $\beta_2=\gcd(a_2, a_4, N)$.
\par
\begin{Prop} \label{Degeneration}
\begin{itemize}
\item[(a)] The degree of $Z_0^1\to P_0^1$ (resp.\ $Z^2_0\to
P^2_0$) is $N/\beta_1$ (resp.\ $N/\beta_2$).
\item[(b)] The genus of $Z^1_0$ (resp.\ $Z^2_0$) is $(N-\gcd(a_1, N)-\gcd(a_3,
  N)-\gcd(a_1+a_3, N))/2\beta_1$ (resp.\ $(N-\gcd(a_2, N)-\gcd(a_4,
  N)-\gcd(a_1+a_3, N))/2\beta_2$). 
\item[(c)] The number of singular points of $\cZZ_0$ is $\gcd(a_1+a_3,
  N)$.
\end{itemize}
\end{Prop}
\par

\section{Veech's $n$-gons revisited} \label{ngonrev}

In this section we realize the $(n,\infty,\infty)$-triangle groups as
the affine groups of a Teichm\"uller curves. This result is due to
Veech, but our method is different. An advantage of our method is that
we obtain the Lyapunov exponents in Corollary \ref{LyapnevenXX} with
almost no extra effort.  The reader may take this section as a
guideline to the more involved next section. In this section the
family of cyclic covers we consider has only one elliptic fixed
point. A $(\ZZ/2\ZZ)$-quotient of this family is shown to be a
Teichm\"uller curve. In the next section there are two elliptic fixed
points and we need a $(\ZZ/2\ZZ)^2$-quotient.  Moreover, common
divisors of $m$ and $n$ in the next section make a fiber product
construction necessary that does not show up here.  
\par 
Let $n=2k\geq
4$ be an even integer and fix a primitive $n$th root of unity
$\zeta_n$.  We specialize the results of \S \ref{covofPP} to the
family $g: \cZZ \to \PP^*$ of curves of genus $n-1$ given by the
equation
$$\cZZ_t:\qquad z^n = x (x-1)^{n-1}(x- t),$$
i.e.\ we consider the case that  $N=n$, $a_1=a_3=1$ and $a_2=a_4=n-1$.
Let $$\varphi(x,y)=(x,\zeta_n y)$$ be a generator of $\Gal(\cZZ/\cPP)$. 
The geometric fibers of $g$ admit an involution covering
$x \mapsto t/x$  $\cPP$ . We choose this
involution to be 
$$ \sigma(x,y) = \left\{\begin{array}{ll}
(\displaystyle{\frac{t}{x}, \frac{t^{2/n} (x-1)(x-t)}{xy}})
&\mbox{ if $k$ is even},\\
(\displaystyle{\frac{t}{x}, \zeta_n \frac{t^{2/n} (x-1)(x-t)}{xy}})
&\mbox{ if $k$ is odd}.
\end{array}\right.
$$
\begin{Lemma} The exponents $a_i$ are chosen such that
\label{Veechlem}
\begin{itemize}
\item[(a)] condition (\ref{Higgseq}) is satisfied for $i=(n+1)/2$,
\item[(b)] the projective monodromy group of the local systems $\LL((n-1)/2))$
and $\LL((n+2)/2)$ is the triangle
group $\Delta(n,\infty,\infty)$.
\end{itemize}
\end{Lemma}
\par {\bf Proof:} Part (a) follows by direct verification. Part (b) is
proved in \cite{CoWo90}.  \hspace*{\fill}$\Box$ \par Let $\pi:C\to \PP^\ast$ be
defined by $s=t^{n/2}$. Then $\sigma$ extends to an automorphism of the
family of curves $g_C: \cZZ_C \to C$. As before, we let $\ol{\pi}:
\ol{C} \to \PP^1$ be the extension of $\pi$ to a smooth
completion. Moreover, the local monodromy matrices of the pullback of
the local systems $\LL(i)$ to $\ol{C}$ are unipotent.  \par We let $f:
\cXX = \cZZ /\langle \sigma \rangle \to C$. Let $\ol{f}: \ol{\cXX} \to
\ol{C}$ be the stable model of $f$.  Our goal is to show that the
fibers $\cXX_c$ of $\ol{f}$ are smooth for all $c\in \widetilde{C} :=
\ol{\pi}^{-1}(\PP^1\sms\{1,\infty\})$. This allows us to apply the
criterion (Proposition \ref{Higgsprop}) for $\widetilde{C}$ to be the
cover of a Teichm\"uller curve. \par
\begin{Thm} \label{Teichneven} Let $g= (n-2)/2$.
The natural map $m: \widetilde{C} \to M_g$ induced by $\ol{f}$
exhibits $\widetilde{C}$ as the unramified cover of a Teichm\"uller
curve.
\end{Thm}
\par {\bf Proof:} We first determine the degeneration of $g_C$ at $c
\in \ol{C}$ with $\ol{\pi}(c) \in \{0,1,\infty\}$. Our assumption on
the local monodromy matrices implies that the fiber $\cZZ_c$ is a
semistable curve, and we may apply Proposition \ref{Degeneration}.
For $\ol{\pi}(c) \in \{1,\infty\}$ the fiber $\cZZ_c$ consists of two
irreducible components, which have genus $0$. The local monodromy
matrices of $\LL(i)$ at $c$ are unipotent and of infinite order for all $i$, as
can be read off from the local exponents.

 Similarly,  the local 
monodromy at $c$ with $\ol{\pi}(c)=0$ is finite. The definition of
$C$ implies therefore that it its trivial. The set $\Su \subset \ol{C}$
(notation of \S \ref{covofPP}) consists exactly of
$\ol{\pi}^{-1}\{1,\infty\}$.  

One checks that $\sigma$ acts on the holomorphic $1$-forms $\omega_i$
(Lemma \ref{dimlem}.(c)) as follows:
\begin{equation} \label{sigmaonone}
\sigma^* \omega_i = (-1)^i d(i) \omega_{n-i}
\quad \text{for} \quad i \neq n/2, \quad\qquad
\sigma^* \omega_{n/2} = - \omega_{n/2},
\end{equation}
 where $d(i) = t^{2i/n-1}$ if $k$ is odd and
$d(i) = t^{2i/n-1} \zeta_n^i$ if $k$ is even.
This implies that the generic fiber of $\cXX$ has genus $n/2-1$.

We claim that $\cXX_c$ is smooth for all $c\in \ol{C}-\Su$. We only
 need to consider $c\in \ol{C}$ such that $\pi(c)=0$.  Proposition
 \ref{Degeneration} implies that the degenerate fiber  $\cZZ_c$
 consists of two components
of genus $n/2-1$. Note that $\sigma$ acts as the permutation
$(0\,\infty)(1\,t)$ on the branch points of $ \cZZ \to \cPP$.  Hence
$\sigma$ interchanges the two components of $\cZZ_c$. We conclude that
the quotient $\cXX_c$ of $\cZZ_c$ by $\rho$ is a smooth curve of genus
$n/2-1$.

Consider the local
system $\MM=\LL((n-2)/2) \oplus \LL((n+2)/2)$ in 
$R^1 (g_C)_* \CC$ on  $C$. It
is invariant under $\sigma$. The part of $\MM$ on which $\sigma$ acts
trivially is a local subsystem $\LL \subset \MM$. This $\LL$ is 
necessarily of rank $2$, since 
$\omega_{(n-2)/2} + d((n-2)/2) \omega_{(n+2)/2}$ 
is $\sigma$-invariant (resp.\ anti-invariant), if $k$ is odd (resp.\ even) 
and $\omega_{(n-2)/2} - d((n-2)/2) \omega_{(n+2)/2}$ is $\sigma$-anti-invariant
(resp.\ invariant) for $k$ odd (resp.\ even). 
This also implies that the compositions
$$ \LL \to \LL((n-2)/2) \oplus \LL((n+2)/2)
\to \LL((n-2)/2) $$
and
$$ \LL \to \LL((n-2)/2) \oplus \LL((n+2)/2) \to \LL((n+2)/2) $$ are
non-trivial. Since the monodromy group, $\Gamma$, of both
$\LL((n-2)/2)$ and $\LL((n+2)/2)$ contains two non-commuting parabolic
elements, we conclude that $\LL((n-2)/2)$ is an irreducible local
system, and hence that
$$ \LL \cong \LL((n-2)/2) \cong \LL((n+2)/2). $$

From Proposition \ref{Higgsprop} and Lemma \ref{Veechlem}.(a) 
we conclude that $\cXX$ is
the universal family over an unramified cover of a Teichm\"uller
curve as claimed.
\hspace*{\fill} $\Box$
\par
 Corollary \ref{LyapnevenXX}  follows from
Proposition\ \ref{lambdalocsyst}:
\par
\begin{Cor} \label{LyapnevenXX}
The VHS of the family $f: \cXX \to C$ decomposes as
$$ R^1 f_* \CC \cong \bigoplus_{j=1}^{(n-2)/2} \LL_j,$$
where $\LL_j$ is a rank $2$ local system isomorphic to $\LL((n-2j)/2)$.
Moreover,  
$$ \lambda(\LL_j) = \frac{k-j}{k-1}.$$
\end{Cor}
\par
Anton Zorich has communicated to the authors that he (with Maxim Kontsevich)
independently calculated these Lyapunov exponents.
\par
\begin{Rem}{\rm
The trace field of $\Delta(n,\infty,\infty)$ is
$K=\QQ(\zeta_n+\zeta_n^{-1})$, hence $r = [K:\QQ] \leq \phi(n/2)$.
Corollary \ref{LyapnevenXX} allows to decomposes the VHS of $\cXX$
completely into rank two pieces. This is much finer than Theorem\
\ref{maxHiggsTeich} that predicts only $r$ pieces of rank two plus
some rest.  }\end{Rem} \par Each fiber $\cZZ_t$ admits an extra
isomorphism, namely
$$\tau(x,y) = \left(\frac{x-t}{x-1}, \;  y \, 
\frac{t-1}{(x-1)^2}
\right)$$
It extends to an automorphism of the family 
$g_{\widetilde{C}} :\cZZ_{\widetilde{C}} \to \widetilde{C}$.
One checks that $\tau$ and $\sigma$ commute. Hence 
$\tau$ descends to an
automorphism of $\cXX$, which we also denote by $\tau$. Let
$p: \cUU = (\ol{\cXX}|_{\widetilde{C}})/\langle \tau \rangle 
\to \widetilde{C}$ the quotient family.
One calculates that
$$ \tau^* \omega_i = (-1)^{i+1} \omega_i.$$ 
From this we deduce that the fibers of $\cXX$ are Veech surfaces
that cover non-trivially Veech surfaces of smaller genus, 
the fibers of the fibers of $p$.
\par
\begin{Thm} \label{Teichnevenquottau}
\begin{itemize}
\item[(a)]
The moduli map $\widetilde{C} \to M_{g(\cUU)}$ of the
family of curves $p: \cUU \to \widetilde{C}$ is an unramified covering of 
Teichm\"uller curve.
Its VHS decomposes as 
$$ R^1 p_* \CC \cong \bigoplus_{j=0}^{t(n)} \LL(1+2j),$$
where $\LL(j)$ is the local system appearing in the VHS of $f$
and $t(n) = (n-6)/4$ if $k$ is odd (resp.\ $t(n) = (n-4)/4$ if $k$ is even).
\item[(b)]
The genus of $\cUU$ is $t(n)+ 1$ and 
$$ \lambda(\LL(1+2j)) = \frac{k-(1+2j)}{k-1}$$
\end{itemize}
\end{Thm}
\par
{\bf Proof:} Both for $k$ odd and $k$ even the generating holomorphic
$1$-form in $\LL(1)$ is $\tau$-invariant. Hence this local system
descends to $\cUU$. The property of being a Teichm\"uller curve now follows
from Proposition \ref{Higgsprop}. The remaining statements are
easily deducted from  Corollary \ref{LyapnevenXX}.
\hspace*{\fill} $\Box$
\par
Let $U$ be a fiber of $\cUU$. We denote 
by $\omega_X \in \Gamma(X, \Omega^1_X)$
(resp.\ $\omega_U \in \Gamma(U, \Omega^1_U)$)
the differential that pulls back to 
$\omega_{(n-2)/2} \pm  d(i) \omega_{(n+2)/2}$ on $\cZZ_c$,
where the sign depends on the parity of $n$ and refer to it as the
{\em generating differential} of the Teichm\"uller curve.
\par
\begin{Cor} \label{isVeechneven}
The Teichm\"uller curve $\cXX$ is the one generated by
the regular $n$-gon studied in \cite{Ve89}.
\end{Cor}
\par {\bf Proof:} Let $c$ be a point of $\ol{C}$ with $\pi(c)=0$.  The
fiber $\cZZ_c$ consists of two components isomorphic to
$$\cXX_0:\qquad y^n=x(x-1)$$ which are interchanged by $\sigma$. The
generating differential $\omega_X$ specializes to the differential
$$\omega_0 =y^{(n-2)/2} \,{\rm d} x/x(x-1)$$ on $\cXX_0$. There is an
obvious isomorphism between the curve $w^n -1 = z^2$ and $\cXX_0$ such
that $\omega_0$ pulls back to the differential ${\rm d} w/z$ considered by
Veech (\cite{Ve89} Theorem 1.1). \hspace*{\fill} $\Box$ \par Actually the
family $\cXX$ is isomorphic (after some base change) to
$$ y^2 = p_t(z) = \prod_{i=1}^n (x - \zeta_n^i - t \zeta_n^{-1}).$$
This was shown by Lochak (\cite{Lo05}, see also \cite{McM04c}). 
\par
The following proposition is shown in \cite{Ve89} Theorem 1.1.
We give an alternative proof in our setting.
\par
\begin{Prop} \label{affgrp1}
The projective affine group of a fiber of $\cXX$ together with 
the generating differential 
contains the $(n,\infty,\infty)$-triangle group.
The same holds for the fibers of $\cUU$.
\end{Prop}
\par
{\bf Proof:} We first consider $\cXX$. We have to show that the moduli map 
$C \to M_g$ given by $\cXX$ factors through $\pi:C \to \PP^*$.
That is, we have to show that
two generic fibers $\cXX_c$ and $\cXX_{\tilde{c}}$ with
$c,\tilde{c} \in C$ such that $\pi(c) = \pi(\tilde{c})$
are isomorphic. Equivalently, we have to show that
for $c,\tilde{c}$ as above there is an isomorphism 
$i_0:\cZZ_c \to \cZZ_{\tilde{c}}$ which is $\sigma$-equivariant.
It suffices to show the existence of $i_0$ after any
base change $\pi: C' \to \PP^*$ such that $\sigma$ is
defined on $\cYY_{C'}$. We may suppose that 
$\ol{\pi}: \ol{C'} \cong \PP_s^1 \to \PP_t^1$ is given by $t = s^{n/2}$.
The hypothesis $\pi(c) = \pi(\tilde{c})$ implies
that $c = \zeta_n^{2j} \tilde{c}$, for some $j$.
It follows that the canonical isomorphism $i: \cZZ_c \to \cZZ_{\tilde{c}}$,
given by $(x,y) \mapsto (x,y)$, satisfies
$$ \sigma \circ i = \varphi^{2j} \circ i \circ \sigma. $$
Hence $i_0 = \varphi^j \circ i $ is the isomorphism we were looking for.
\par
The proof for the family $\cUU$ is similar.
\hspace*{\fill} $\Box$
\par
We record for completeness:
\par
\begin{Cor}
All $(n,\infty,\infty)$-triangle groups for $n \geq 4$ arise as projective
affine groups.
\end{Cor} 
\par
\begin{Rem} \label{oddVeech}
{\rm For $n$ odd the same construction works with $N$ and
$a_i$ chosen as above. The local exponents of $(\LL(i),\omega_i)$ 
at $t=0$ are then $1-2i/n$. The local system $\LL(i_0)$ becomes 
maximal Higgs for  $i_0=(n+1)/2$, after a base 
change $\pi: C \to \PP^*$ 
whose extension to $\ol{C} \to \PP^1$ is branched of order $n$
at $0$. The quotient family $f: \cXX = \cZZ/\langle \sigma \rangle \to C$
may be constructed in the same way as above. Its moduli map yields as above
a Teichm\"uller curve  $\widetilde{C }\to M_g$ where $g=(n-1)/2$. 
The corresponding translation 
surfaces are again the ones studied in \cite{Ve89}. Veech also
determines that the affine group is not $\Delta(n,\infty,\infty)$
but the bigger group $\Delta(2,n,\infty)$, containing $\Delta(n,\infty,\infty)$
with index two. We obtain the same family of curves also
as a special case of the construction in \S \ref{realize}, 
by putting $m=2$. For this family we calculate, using 
Proposition\ \ref{lambdalocsyst}, that
$$ \lambda(\LL(i)) = \frac{2i}{n-1}, \qquad i=1,\ldots,(n-1)/2.$$
} \end{Rem}

\section{Realization of $\Delta(m,n,\infty)$ as projective affine group} 
\label{realize}

Let $m,n>1$ be integers with $mn \geq 6$. We let
\[
\sigma_1=\frac{nm+m-n}{2mn},\quad \sigma_2=\frac{nm-m+n}{2mn},\quad
\sigma_3=\frac{nm+m+n}{2mn},\quad \sigma_4=\frac{nm-m-n}{2mn}.
\]
and we let $N$ be the least common denominator of these fractions.
We let $a_i = N \sigma_i$ and consider the
family of curves $g: \cZZ \to \PP^*$ given by
$$ \cZZ_t: \qquad y^N = x^{a_1}(x-1)^{a_2}(x-t)^{a_3}. $$
The family $g$ cyclically covers the constant family $\cPP \cong
\PP^1_x \times \PP^* \to \PP^*(=\PP^1_t-\{0,1,\infty\})$.
\par
The plan of this section is as follows. We construct a cover $\cYY
\to \cZZ$ such that the involutions 
 \begin{equation}\label{eq1}\begin{split}
\sigma(x)&=
  \left(t(x-1)/(x-t)\right),\\
  \tau(x)&=\left(t/x\right)
\end{split}
\end{equation}
of $\cPP \to \PP^*$ lift to involutions of the family $\cYY_C \to C$
obtained from $\cYY \to \PP^*$ by a suitable unramified base change 
$\pi: C \to \PP^*$.  We denote these lifts again by $\sigma$ and $\tau$.
If $m$ and $n$ are relatively prime then in fact $\cYY$ equals $\cZZ$.
\par
\begin{Rem}{\rm
The exponents $a_i$ are chosen such that the local system $\LL_{\chi}$
has as projective monodromy group the 
triangle group $\Delta(m,n,\infty)$, see again e.g.\ \cite{CoWo90}.
We modify the lifts $\tau$ and $\sigma$ by appropriate powers of 
a generator of 
$\Aut(\cZZ/\cPP)$ such that the group $H=\langle \tau, \sigma \rangle$ 
is still isomorphic
to $(\ZZ/2)^2$ and such that $\sigma$ and $\tau$ and $\sigma\tau =: \rho$
have `as many fixed points as possible'.
}\end{Rem}
\par
We  consider the quotient
family $f: \cXX = \cYY/H \to C$. Its stable model 
$\ol{f}: \ol{\cXX} \to \ol{C}$ has smooth fibers over 
$\widetilde{C} = \ol{\pi}^{-1}(\PP^1\sms\{\infty\})$, 
where $\ol{\pi}: \ol{C} \to \PP^1$ extends $\pi$.
\par
Together with an analysis of the action of $H$ on differentials
we can apply Proposition \ref{Higgsprop} to produce Teichm\"uller curves.
\par 
\begin{Thm} \label{hisTeich}
Via the natural map $m: \widetilde{C} \to M_g$ induced from $\ol{f}$
the curve $\widetilde{C}$ is an unramified cover of a Teichm\"uller curve.
The genus $g$ is given in Corollary\ \ref{genuscor}.
\end{Thm}
\par As corollaries to this result we calculate the precise VHS of
$f$ and the projective affine group of the translation surfaces
corresponding to $f$.  In \S \ref{Wardsec} we show that for $m=3$
we rediscover Ward's Teichm\"uller curves (\cite{Wa98}).  \par
\begin{Rem}\label{simplerem}{\rm 
The notation in the proof of Theorem \ref{hisTeich} is rather
complicated, due to the necessary case distinction. We advise the
reader to restrict to the case that $m$ and $n$ are odd and relatively
prime on a first reading. This considerably simplified the notation,
but all main features of the proof are already visible. In this case
$\cYY=\cZZ$, and  $m=m'$, $n=n'$,
$\gamma=\beta=1$, and $N=\bar{N}=\hat{N}$.  }\end{Rem} \par We start
with some more notation. We write $Z$ (resp.\ $P$, $X$, $Y$) for the
geometric generic fiber of $\cZZ$ (resp.\ $\cPP$, $\cXX$, $\cYY$).  We
choose a primitive $N$th root of unity $\zeta_N\in \CC$ and define the
automorphism $\varphi_1 \in \Aut(\cYY/\cPP)$ by
$$ \varphi_1(x,y) = (x, \zeta_N y). $$
\par
We need to determine the least common denominator $N$ of the $\sigma_i$, 
$i=1,\ldots,4$, precisely. Let $m = 2^\mu m'$, $n=2^\nu n'$ with $m',n'$ odd. 
We may suppose that $\mu \geq \nu$. Define
$$ \gamma_1 = \gcd(2mn,mn+m-n), \quad \gamma_2 = \gcd(2mn,mn+m+n),
\quad \gamma = \gcd(m,n) $$ and write $\gamma = 2^\nu \gamma'$.  We
distinguish four cases and determine $N =
2mn/\gcd(\gamma_1,\gamma_2)$, accordingly.
$$\begin{array}{lcll}
\text{Case O: odd} & \mu = \nu = 0, & N = 2mn/\gamma, & 
\widehat{N} = N/\gamma= 
2^\delta m'n'/\gamma'^2,\\
\text{Case OE: $m$ odd, $n$ even} & \mu > \nu =0 , & N = 2mn/\gamma,
& \widehat{N} = N/\gamma = 2^\delta m'n'/\gamma'^2, \\
\text{Case DE: different $2$-val., even} & \mu > \nu >0, & N = 2mn/\gamma,
& \widehat{N} = 2N/\gamma = 2^\delta m'n'/\gamma'^2, \\
\text{Case S: same $2$-valuation, even} & \mu = \nu \neq 0, & N=mn/\gamma, &
\widehat{N} = N/\gamma= mn/\gamma^2. \\
\end{array}
$$ It is useful to keep in mind that $\gamma =
\gcd(\gamma_1,\gamma_2)$, except in case $S$ where $2\gamma =
\gcd(\gamma_1,\gamma_2)$. We let $\delta := 0$ in case S, and $\delta
:= \min\{\mu-\nu+2, \mu+1\}$, otherwise.  \par Our first goal is to
determine the maximal intermediate covering of $Z \to P$ to which
$\tau$ lifts. This motivates the definition of $\widehat{N}$.  Let
$0<\bar{\alpha}< \widehat{N}$ be the integer satisfying
\[
\bar{\alpha}\equiv 1\bmod{m'/\gamma'}, \quad \bar{\alpha}\equiv 
-1\bmod{n'/\gamma'}, \quad \bar{\alpha} \equiv \left\{ 
\begin{array}{lll} 1 & \bmod \,2^\delta & \text{cases O, OE, S,} \\
\displaystyle{\frac{n' + 2^{\mu-\nu}m'}{n' - 2^{\mu-\nu}m'}} 
&\bmod \, 2^\delta & \text{case DE}. \end{array}\right. 
\]
For convenience, we  lift $\bar{\alpha}$ to an element $\alpha$ in
$\ZZ/N\ZZ$ such that $\alpha^2=1$. 
\par
Recall that for a rational number $\sigma$, we write $\sigma(i):=\langle
i\sigma\rangle$ (the fractional part). Similarly, for an integer $a$ we write
$a(i)= a(i;\nu)=\nu\langle ia/\nu\rangle$, where $\nu$ is mostly clear from
the context. For each integer $0<i<N$ which is prime to $N$, 
we write 
\[
z(i)=\frac{z^i}{x^{[i\sigma_1]}(x-1)^{[i\sigma_2]}(x-t)^{[i\sigma_3]}},
\quad \text{hence} \quad
z(i)^N=x^{a_1(i)}(x-1)^{a_2(i)}(x-t)^{a_3(i)}.
\]
\par
\begin{Lemma}\label{autlem1}
\begin{itemize}
\item[(a)]
In the cases O, OE and DE the covering $Z \to P$ has ramification order 
$\gamma N/\gamma_1$ (resp.\ $\gamma N/\gamma_2$)
in  points of $Z$ over $x=0, 1$ (resp.\ $x=t, \infty$). 
In case $S$ the ramification orders are $\gamma N/2\gamma_1$  
(resp.\ $\gamma N/2\gamma_2$). Therefore 
$$ 
g(Z)= \left\{ \begin{array}{ll} 
1+N-{\gamma_1+\gamma_2}/{2\gamma} & \text{case S}, \\
1+N-{\gamma_1+\gamma_2}/{\gamma}  & (\text{other cases}).
\end{array} \right. 
$$
\item[(b)]
The automorphism $\sigma$ of $P$ lifts to an automorphism  $\sigma$ of $Z$ of
order $2$.
\item[(c)] The automorphism $\tau$ of $P$ lifts to an automorphism 
$\tau$ of order $2$ of $\widehat{Z}:=Z/ \langle\varphi_1^{\widehat{N}}
\rangle$. 
Moreover, we may choose the lifts such that $\sigma, \tau$ commute as 
elements of $\Aut(\widehat{Z})$.
\item[(d)] We may choose the lifts $\sigma, \tau$ such that,
moreover, $\tau$ has $4m/\gamma$ fixed points (resp.\
$2m/\gamma$ in case S) and such that  $\rho := \sigma \tau$ has
$4n/\gamma$  fixed points on $\widehat{Z}$ (resp.\ $2n/\gamma$ in case S).
\item[(e)] With $\sigma$ and $\tau$ chosen as in (d) the 
automorphism $\sigma$ has no ($2$ in case S) fixed points both
on $Z$ and on $\widehat{Z}$.
\end{itemize}
\end{Lemma}
\par
{\bf Proof:} The statements in (a) are immediate from the
definitions. For (b) and (c) we choose once and for all elements $t^{1/n},
 (t-1)^{1/m}\in\overline{\CC(t)}$.  Define  
\begin{equation}\label{cdeq}
c=(t-1)^{\sigma_2+\sigma_3}, \qquad
d=t^{\sigma_1+\sigma_3}.
\end{equation}
Then
\[
\sigma(z)=c d \frac{x(x-1)}{z (x-t)}= c d \frac{z(-1)}{(x-t)^2}
\]
defines a lift of $\sigma$ to $Z$, since
$\sigma_1+\sigma_2=\sigma_3+\sigma_4=1$. Moreover, this lift has order $2$. We
denote it again by $\sigma$. The quotient curve $\widehat{Z}$ is defined by 
the equation
\[
\bar{z}^{\widehat{N}}=x^{\bar{a}_1}(x-1)^{\bar{a}_2}(x-t)^{\bar{a}_3},
\]
where $\bar{a}_i$ denotes $a_i \bmod \widehat{N}$. 
One computes that $\alpha$ satisfies:
\begin{equation}\label{aitaueq}
(\bar{a}_1(\alpha), \bar{a}_2(\alpha), \bar{a}_3(\alpha), \bar{a}_4(\alpha))=
(\bar{a}_4, \bar{a}_3, \bar{a}_2, \bar{a}_1).
\end{equation}
This implies that 
\[
\tau(\bar{z})=d^\gamma \frac{\bar{z}(\alpha)}{x^{2\gamma}}
\]
defines a lift of $\tau$ to $\widehat{Z}$ which has order $2$. It is easy to
check that $\tau$ commutes with the image of $\sigma$ on $\widehat{Z}$. This
proves (b). Furthermore, one checks that $\sigma$ is an involution and that
$$ \tau \varphi_1 \tau = \varphi_1^{\alpha} \in \Aut(\widehat{Z})
\quad \text{and} \quad \sigma \varphi_1 \sigma = \sigma^{-1} \in
\Aut(Z). $$ This proves (c).
 \par We start with the proof of (d).  Let $x_1 = \sqrt{t}$
be one of the fixed points of $\tau$ on $P$ and let $R$ be a point in
the fiber of $\widehat{Z} \to P$ over $x_1$.  We may describe the
whole fiber by $R_a := \varphi_1^a R$ for
$a=0,\ldots,\widehat{N}-1$. Suppose that $\tau R = R_{a_0}$, hence
$\tau R_a = R_{a_0 + \alpha a}$. Since $\tau$ is an involution, $a_0$
satisfies necessarily $a_0 \equiv 0 \bmod m'/\gamma'$ and $2a_0 \equiv
0 \mod 2^\delta$. Furthermore, $R_a$ is a fixed point of $\tau$ if and
only if
\begin{equation} \label{condtau}
a_0 \equiv 2 a \bmod n'/\gamma' \quad \text{and} \quad 
a_0 \equiv 2^{\mu-\nu+1} a \frac{-m'}{n' - 2^{\mu-\nu}m'} \mod 2^\delta.
\end{equation}
Hence if $\tau$ has a fixed point in this fiber it has precisely
$2^{(\mu-\nu+1)} m'/\gamma'$ 
fixed points in this fiber ($m'/\gamma' = m/\gamma$ in case S).
Since $\tau$ and $\sigma$ commute, $\sigma$ bijectively
maps fixed points of $\tau$ over $x_1$ to fixed points of $\tau$ over 
$x_2 = -\sqrt{t}$. Hence, if $\tau$ has a fixed point, then the number of
fixed points is as stated in (d).

Similarly, let $x_3 = 1+\sqrt{1+t}$ be one the fixed points of $\rho$
on $P$ and let $S$ be a point in the fiber over $x_3$. Write
$S_b=\varphi_1^b S$ for the whole fiber.  Write $\rho S = S_{b_0}$. As
above we deduce that $b_0\equiv 0 \bmod m'/\gamma'$ and $2^{\mu-\nu+1}
b_0 \equiv 0 \bmod 2^\delta$. Then $S_b$ is a fixed point of $\rho$ if
\begin{equation} \label{condrho}
b_0 \equiv 2 b \bmod m'/\gamma' \quad \text{and} \quad
b_0 \equiv 2 b \frac{n'}{n' - 2^{\mu-\nu}m'} \bmod 2^\delta.
\end{equation}
Analogously to the argument for $\tau$, one checks that if $\rho$
has a fixed point then it has as many fixed points as claimed in (d).

Note that we may replace $\sigma$ by $\varphi^i \sigma$
and $\tau$ by $\varphi^j \tau$ without changing the orders of these
elements and such that they still commute if the following conditions
are satisfied: 
\begin{equation} \label{condKleinV}
j \equiv 0 \bmod m'/\gamma', \quad  \quad j \equiv i \bmod n'/\gamma'
\quad \text{and} \quad 2j \equiv 2^{\mu-\nu+1} i \bmod 2^\delta.
\end{equation}
\par
The only obstruction for $\tau$ and $\rho$ to have fixed points consists
in the condition modulo $2^\delta$. We check in each case that
we can modify  $\tau$ and $\rho$ respecting (\ref{condKleinV}) such
that this obstruction vanishes.
\par
In case S there is nothing to do, since $\delta =0$. In case O we might 
have to change the parity of $a_0$ and $b_0$ or both, since $\delta=1$. This
is possible since  (\ref{condKleinV}) imposes no parity condition
in this case: we replace $\sigma$ by $\varphi^i \sigma$ and
$\tau$ by $\varphi^j \tau$ such that $j \equiv a_0 \bmod 2$ and
$i+j \equiv b_0 \bmod 2$. In case OE the conditions for $\tau$ to have
fixed points are satisfied. We might have to change the parity of
$b_0$ which can be achieved since (\ref{condKleinV}) imposes no
parity conditions on $i$ in this case. In case DE we can
solve equations (\ref{condtau}) (resp.\ (\ref{condrho})) for $a$ (resp.\ $b$)
using the conditions imposed on $a_0$ and $b_0$ from the assumptions that
 $\tau$ and $\rho$ are involutions. This proves (d).
\par
For (e) we check with the same argument as above that $\sigma$
has $0$ or $4$ (resp.\ $0$ or $2$ in case S) fixed points.
Checking case by case one finds that $\widehat{Z} \to P$ is totally ramified 
over $\{0,1,t,\infty\}$. Hence $g(\widehat{Z}) = \widehat{N} -1$. The 
Riemann--Hurwitz formula implies that $\sigma$ does not have   fixed points on 
$\widehat{Z}$, hence also not on $Z$ in 
case O, D and DE. The number of fixed points of $\sigma$ in case S 
may be checked directly by counting
 fixed points of $\tau$ on $\widehat{Z}$.
\hspace*{\fill} $\Box$
\par
Let $Z^\tau$ be the conjugate of $Z$ under $\tau$. 
Define $Y$ as the normalization of
$Z\times_{\widehat{Z}} Z^\tau$.  As remarked above, the definition of 
$\widehat{N}$ implies that $\widehat{Z}\to P$ is the largest subcover of 
$Z\to P$ such that $\tau$ 
lifts to $\widehat{Z}$. In other words, $Y\to \widehat{P} := P/\langle
\sigma, \tau \rangle$ is the Galois closure of $Z\to \widehat{P}$. This implies
that $Y$ is  connected. I.e., the particular choice of $\widehat{N}$ is
used precisely to guarantee that the Veech surfaces constructed in
Theorem \ref{hisTeich} are connected. 
\par
By construction, $\sigma$ lifts to $Z$ 
acting on both $Z$ and $Z^\tau$ and $\tau$ lifts to $Z$ by exchanging
the two factors of the fiber product. These two involutions commute
and $\rho:=\sigma\tau$ also has order $2$. 
We have defined the following coverings. The labels indicate the Galois group
of the morphism with the notation introduced in the following lemma. 
\[
\xymatrix{
&Y \ar[dl]_{\langle \psi_2 \rangle}
\ar[dr]^{\langle \psi_1^{\widehat{N}} \psi_2^{-1}  \rangle}&\\
Z\ar[dr]_{\langle\varphi^{\widehat{N}}_1\rangle}&&Z^\tau\ar[dl]^{\langle 
\varphi_2^{\widehat{N}}\rangle}\\
&\widehat{Z}\ar[d]^{\langle \varphi_1 \bmod \widehat{N} \rangle 
= \langle \varphi_2 \bmod \widehat{N} \rangle}&\\
&P \ar[d]^{\langle \sigma,\tau \rangle}\\
&\widehat{P}
}
\]
\begin{Lemma}\label{autlem2}
\begin{itemize}
\item[(a)] We may choose a generator $\varphi_2$ of $\Aut(Z^\tau/P)$
such that the Galois group, $G_0$, of $Y/P$ is
$$ G_0 \cong \{(\varphi_1^i,\varphi_2^j), \,\,i,j, 
\in \ZZ/N\ZZ, \, i \equiv j \bmod \widehat{N}  \} 
\subset \langle \varphi_1 \rangle \times
\langle \varphi_2 \rangle \cong (\ZZ/N\ZZ)^2. $$
We fix generators $\psi_1 = (\varphi_1, \varphi_2)$ and
$\psi_2 = (0, \varphi_2^{\widehat{N}})$ of $G_0$. The Galois group, $G$,
of the covering
 $Y/\widehat{P}$ is generated by $\psi_1, \psi_2, \sigma, \tau$, satisfying
\[
\psi_1^{N}=\psi_2^{\beta} =\sigma^2=\tau^2=1,\qquad
 [\psi_1, \psi_2]=[\sigma, \tau]=1,
\]
\[
\sigma\psi_i \sigma=\psi_i^{-1}\quad (i=1,2), \qquad 
\tau\psi_1 \tau=\psi_1^{\alpha}, \qquad 
\tau\psi_2 \tau=\psi_1^{\alpha N} \psi_2^{-\alpha} 
(= (\varphi_1^{\alpha N}, 0)).
\]
\item[(b)] The genus of $Y$ is $g(Y) = 1+N\beta-2\beta$, 
where $\beta = \gamma/2$ in case DE and $\beta = \gamma$ in the
other cases.
\item[(c)] The number of fixed points of $\tau$ on $Y$ is $4m\beta/\gamma$ 
(resp.\ $2m$ in
case S).
\item[(d)] The number of fixed points of $\rho$ on $Y$ is $4n\beta/\gamma$ 
(resp.\ $2n$ in case S).
\item[(e)] The involution $\sigma$ has no fixed points on $Y$.
\end{itemize}
\end{Lemma}
\par
{\bf Proof:} The presentation in (a) follows from the above construction.
To prove (b), we remark that $Z^\tau$ is given by the equation
\[
\tilde{z}^{N}=x^{a_4}(x-1)^{a_3}(x-t)^{a_2},
\]
compare to (\ref{aitaueq}).
Recall that $\widehat{Z} \to P$ is totally ramified over
$\{0,1,t,\infty\}$. Hence at each of the $\gamma_1/\gamma$ points (resp.\
$\gamma_1/2\gamma$ in case $S$) over $0$ and $1$ in $Z$
the map $Z\to \widehat{Z}$ is branched of order $\gamma^2/\gamma_1$ 
(resp.\ $2\gamma^2/\gamma_1$ in case S and $\gamma^2/2\gamma_1$ in
case DE). 
The other covering $Z^\tau\to \widehat{Z}$ is branched at the corresponding
$\gamma_1/\gamma$ (resp.\ $\gamma_1/2\gamma$ in case $S$) 
points of order $\gamma^2/\gamma_2$ (resp.\ $2\gamma^2/\gamma_2$ 
in case S and and $\gamma^2/2\gamma_2$ in
case DE). Over $t$ and $\infty$ instead of $0$ and $1$ 
the roles of $\gamma_1$ and $\gamma_2$ are interchanged.
\newline
It follows from Abhyankar's Lemma that $Y\to \widehat{Z}$ is ramified 
in all cases at each point over $0,1,t,\infty$ of order
$\beta$. Hence these fibers of $Y \to P$ consist of $\beta$ 
points in each case.
\par
For (c), (d) and (e) note that $Z \to P$ is unramified over
the fixed points of $\tau$, $\sigma$ and $\rho$. Hence $Y$ is
indeed the fiber product in neighborhoods of these points.
Since $\tau$ interchanges the two factors, exactly $\beta$ of the
$\beta^2$ preimages in $Y$ of a fixed point of $\tau$ on $Z$ 
will be fixed by the lift of $\tau$ to $Y$. This completes 
the proof of (c). 
\par
For (d) note that $\id \times \sigma: Z \times_{\widehat{Z}} Z^{\tau}
\to  Z \times_{\widehat{Z}} Z^{\sigma}$ is an isomorphism and we may
now argue as in (c).
\par
If $\sigma$ has a fixed point on $Y$ it has a fixed point on $Z$. 
This implies (e) for cases O, OE and DE. In case S we argue as in
the proof of Lemma \ref{autlem1}, and conclude that $\sigma$ has $0$ or two
 fixed points
in $Y$ above each fixed point in $\widehat{Z}$. We deduce the claim from 
the Riemann--Hurwitz formula applied to $Y \to Y/H$ .
\hspace*{\fill} $\Box$
\par
\begin{Cor} \label{genuscor}
The genus of $X=Y/H$ is $g(X)= (mn-m-n - \gamma)\beta/ 2\gamma+1$ 
in case O, OE and D and $g(X) = (mn-m-n-\gamma)/4 +1$ in case S.
\end{Cor}
\par
\begin{Not}\label{pinot}
  {\rm Until now we have been working on the geometric generic fiber of $g:
    \cYY \to \PP^*$ etc.  Let $\pi: C \to \PP^*$ be the unramified cover
    obtained by adjoining the elements $c,d$ defined in (\ref{cdeq}) to
    $\CC(t)$. Then $H=\langle\sigma, \tau\rangle$ is a subgroup of
    $\Aut(\cYY_C)$.  Passing to a further unramified cover, if necessary, we
    may suppose that the VHS of the pullback family $h_C: \cYY_C \to C$ is
    unipotent.  We write $\ol{\pi}:\ol{C}\to \PP^1_t$ for the corresponding
    (branched) cover of complete curves. Then $h_C$ extends to a family $h_C:
    \ol{\cYY}_C \to \ol{C}$ of stable curves over this base curve.  }\end{Not}
  \par The following lemma describes the
action of $H$ on the degenerate fibers of $h_C$.
\begin{Lemma} \label{singfibsmooth}
Let $c \in \ol{C}$ be a point with $\pi(c) \in \{0,1\}$. The quotient
$\cXX_c := (\cYY_C)_c/ H$ is smooth and
$$ 
 g(\cXX_c) = \left\{
\begin{array}{ll}
(mn-m-n - \gamma)\beta/ 2\gamma+1 & 
\text{cases O, OE and DE}, \\
(mn-m-n-\gamma)/4 +1  & \text{case S}.\\
\end{array} \right. $$
\end{Lemma}
\par
{\bf Proof:}
Choose $c\in \pi^{-1}(0)$. The case that $c\in \pi^{-1}(1)$ is similar, and
left to the reader.  By Proposition\ \ref{Degeneration}
the fiber $(\cZZ_C)_c$ consists of two irreducible
components which we call $Z_0^1$ and $Z_0^2$; we make the convention that
the fixed points $x=0, t$ of $\varphi_1$ on $\cZZ_C$ specialize to $Z_0^1$.
Choosing suitable coordinates, the curve $Z_0^1$ is given by
\begin{equation}\label{Z01eq}
z_0^{{N}}=x_0^{a_1}(x_0-1)^{a_3}.
\end{equation}
The components $Z_0^1$ and $Z_0^2$ intersect in $2m/\gamma$ points 
(resp.\ $m/\gamma$ in case S). We write
$P_0^j$ for the quotient of $Z_0^j$ by $\langle \varphi_1 \rangle \cong \ZZ/N$.
\par 
We claim that the fiber $(\cYY_C)_c$ consists of $2$ irreducible 
components $Y_0^1, Y_0^2$, as well. Let ${\mathcal N}$ 
be the normalization of the
fiber product $(\cZZ_C)_c \times_{(\bar{\cZZ}_C)_c}(\cZZ_C)^\tau_c$. 
By Abhyankar's Lemma,  ${\mathcal N} \to (\cZZ_C)_c$ is \'etale at the
preimages of the intersection point of the two components of $(\cPP_C)_c$.
Hence $N$ consists of two curves: the fiber products over $Z^j_0/\langle
\varphi_1^{\widehat{N}}\rangle$ of $Z_0^j$ with its $\tau$-conjugate, for
$j=1,2$. These two curves intersect transversally in $2m\beta/\gamma$ points.
This implies that ${\mathcal N}$ is a stable curve and indeed the fiber $(\cYY_C)_c$.
\par
One computes that $g(Y_0^j)=1 + mn - m\beta/\gamma - \beta$ 
in cases O, OE and DE
and $g(Y_0^j)= mn -m/2+1-\gamma$ in case $S$. Since $\rho$ acts on the
points $\{0,1,t,\infty\}$ as the permutation $(0\,t)(1\,\infty)$ we
conclude that $\rho$ fixes the components $Y_0^j$ while $\sigma$
and $\tau$ interchange them. Clearly, for a coordinate $x_0$ as in 
(\ref{Z01eq}) we have that $\rho(x_0) = 1-x_0$, i.e.\ $\rho$ fixes
the points $1/2$. This is a specialization of one of the two
fixed points $1\pm\sqrt{1-t} \in P$. Since by Lemma \ref{autlem2}
the automorphism $\rho$ fixes $2n$ ($n$ in case $S$) points in
$Y$ above each of these points of $P$ it follows that $\rho$ fixes
$2n$ (resp.\ $n$) points of $Y_0^j$ with $x_0=1/2$. It remains to
compute the number, $r_\infty$, of fixed points of $\rho$ over $x_0=\infty$.
\par
Suppose we are not in case S. Then by the Riemann--Hurwitz formula
$$g(\cXX_c) = g(Y_0^j/\langle \rho \rangle) = 
(mn-m-n -\gamma)\beta/2\gamma +1 - r_\infty/4.$$
Applying the Riemann--Hurwitz formula to the
quotient map $Z_0^j \to Z_0^j/\langle \rho \rangle$, we conclude that 
$r_\infty \equiv 0 \bmod 4$. Represent
the fiber of $Z_0^j$ over $\infty$ as $\varphi_1^b R$, for 
$b=1,\ldots,2m/\gamma$. As in the proof of Lemma \ref{autlem1}, we conclude 
that $r_\infty$ equals
zero or two. It follows that $r_\infty =0$. 
\newline
In case S we have
$$g(\cXX_c) = (mn-m-n-\gamma)/4 +1 - r_\infty/4.$$
and we conclude as above that $r_\infty =0$.
\newline
Genus comparison shows that the fiber $(\cZZ_C)_c$ is smooth.
\hspace*{\fill} $\Box$
\par
{\bf Proof of Theorem\ \ref{hisTeich}:}
We have shown in  Lemma \ref{singfibsmooth} that
$\ol{\cXX}_c$ is smooth for $c\not\in\Su=\pi^{-1}(\infty)$. We have to show 
that the VHS of $f:\cXX\to C$
contains a local subsystem of rank $2$ which is maximal Higgs.
\newline
We decompose the VHS of $g$ into the characters
$$ \chi(i,j): \left\{ \begin{array}{lcl} G_0 & \to &\CC \\
\psi_1 & \mapsto & \zeta_N^i \\
\psi_2 & \mapsto & (\zeta_N^{\widehat{N}})^j. \\
\end{array}
\right. $$ We let $\LL(i,j) \subset R^1 h_*\CC$ be the local system on
which $G$ acts via $\chi(i,j)$. Local systems with $j=0$ arise as
pullbacks from $\cZZ$. By Lemma \ref{dimlem} the local systems
$\LL(i,0)$ are of rank two if $i$ does not divide $N$. Using the
presentation of $G$ one checks that $\sigma^* \LL(i,j) = \LL(-i,-j)$
and $\tau^* \LL(i,j) = \LL(-\alpha i, \alpha(i-j))$.  
\par The local
exponents of $(\LL(1,0), \omega_1)$ at $t=0$ (resp.\ $t=1$) are $(0,
1/n)$ (resp.\ $(0,1/m)$). Therefore, the definition of $
\ol{\pi}:\ol{C}\to \PP^1_t$ (Notation \ref{pinot}) implies that  condition 
(\ref{Higgseq}) is satisfied for $\LL(1,0)$.

Consider the local system 
$$\MM := \LL(1,0) \oplus \LL(-1,0)
\oplus \LL(-\alpha, \alpha) \oplus \LL(\alpha,-\alpha) $$
on $\cZZ_C$.
Since $H$ permutes the $4$ factors of $\MM$ transitively, we conclude 
that for each character $\xi$ of $H$ there is a rank two local
subsystem of $\MM$ on which $H$ acts via $\xi$. Moreover the projection
of the subsystem $\LL := \MM^H$ to each summand is 
non-trivial. Since the $4$ summands of $\MM$ are irreducible by
construction, this implies that
$$ \LL \cong   \LL(1,0) \cong \LL(-1,0)  \cong  
\LL(-\alpha, \alpha) \cong \LL(\alpha,-\alpha).$$
Hence $\LL$ descends to $\cXX$ and is maximal Higgs with respect to $\Su$.
Proposition \ref{Higgsprop} implies that 
the extension of $f$ to $\pi^{-1}(\PP^1\sms\{\infty\})$ is the pullback of 
universal family of curves to an unramified
cover of a Teichm\"uller curve.
\hspace*{\fill} $\Box$
\par
The proof of Theorem \ref{hisTeich} contains more information 
on the VHS of $f$ and on the
Lyapunov exponents $\lambda(\LL_i)$.
We work out the details in the most transparent case that $m$, $n$ are odd 
integers which are relatively prime. The interested reader can easily work out 
the Lyapunov 
exponents in the remaining cases, too. In this case the curves $\cZZ$ and 
$\cYY$
coincide (Remark \ref{simplerem}) and the local system $\LL(i,j)$ is $\LL(i)$
in the
notation of Lemma \ref{dimlem}.
\par
We deduce from the arguments of the  proof of Theorem \ref{hisTeich} that, for 
each 
$i$ not divisible by $m$ or $n$, there is an $H$-invariant local 
system $\LL_i$ with
$$ \LL_i \cong  \LL(i) \cong 
\LL(\alpha i)  \cong  \LL(-\alpha i)  
\cong  \LL(-i).$$ 
Since those $i$ fall into $(m-1)(n-1)/2$ orbits under $\langle
\pm 1, \pm \alpha\rangle$, we have the complete description of 
the VHS of $h$. Let $c_j(i)=\sigma_j(i)+\sigma_3(i)-1$.
\par
\begin{Cor} \label{splitVHSh}
Let $m$ and $n$  be odd integers which are relatively prime.
\begin{itemize}
\item[(a)]
The VHS of $f$ splits as
$$ R^1 f_* \CC \cong \bigoplus_{j \in J} \LL(j),$$
where $\LL(j)$ is an irreducible rank two local system 
and $j$ runs through a set of representatives of  
$$ J = \{0<i<N,  m\nmid i, n\nmid i \}/\sim, \quad
\text{where} \quad i \sim -i \sim \alpha i \sim -\alpha i. $$ 
\item[(b)] The Lyapunov exponents are 
$$ \lambda(\LL(i)) = \frac{mn - e_1(i)m - e_2(i)m}{mn-m-n}, \quad 
\text{where} \quad e_1(i)= n|c_1(i)| \quad\text{and}\quad
e_2(i) = m|c_2(i)|.$$
\end{itemize}
\end{Cor}
\par
{\bf Proof:} This follows directly from Proposition \ref{lambdalocsyst}.
\hspace*{\fill} $\Box$ 
\par
\par
\begin{exa} \label{lyapex} {\rm
We calculate the Lyapunov exponents explicitly for $m=3$ and $n=5$.
Then $N=2nm=30$ and hence $\alpha=19$. We need to
calculate the $\lambda(\LL(i))$ only up to the relation `$\sim$'
and hence expect at most $4$ different values. One checks:
\[
\lambda(\LL(i))=\left\{\begin{array}{ll}
7/7& \mbox{ if }\quad i\sim 1,\\
4/7&  \mbox{ if }\quad i\sim 2,\\
2/7& \mbox{ if }\quad i\sim 4,\\
1/7& \mbox{ if }\quad i\sim 7.
\end{array}\right.
\]
In particular, we see that the $\lambda(\LL(i))$ do in general not form an 
arithmetic
progression as one might have guessed from studying Veech's $n$-gons. 
}\end{exa}
\par
\begin{Rem}\label{splitcompletely}{\rm
Note that  $K:= \QQ(\cos(\pi/n),\cos(\pi/m))$ is the trace field of the 
$\Delta(m,n,\infty)$-triangle group. Hence $r = [K:\QQ] 
\leq \phi(mn)/4 \leq (m-1)(n-1)/4$. Here again the decomposition of the 
VHS is finer than predicted by Theorem\ \ref{maxHiggsTeich},
compare the remark after Corollary \ref{LyapnevenXX}.
}\end{Rem}
\par
Let $X$ be any fiber of $f$. We denote 
by $\omega_X \in \Gamma(X,\Omega^1_{X})$ 
a generating differential, i.e.\ a holomorphic differential
that generates $(1,0)$-part of the maximal Higgs local system 
when restricted to the fiber $X$. This condition 
determines $\omega_X$ uniquely up to scalar multiples.
\par
\begin{Thm} \label{mnwhichtriangle}
The projective affine group of the translation surface $(X,\omega_X)$
is \begin{itemize}
\item[(a)] the $(m,n,\infty)$-triangle group, if $m\neq n$.
\item[(b)] the $(m,m,\infty)$-triangle group or the
$(2,m,\infty)$-triangle group, if $m=n$. 
\end{itemize}
\end{Thm}
\par
{\bf Proof:} We first show that the triangle group $\Delta(m,n,\infty)$
is contained in the projective affine group of $(X,\omega_X)$. 
As in the proof of Proposition\ 
\ref{affgrp1}, we take two fibers $\cYY_{c}$ and $\cYY_{\tilde{c}}$ with
$\pi(c) = \pi(\tilde{c})$. We need to show the existence
of an isomorphism $i_0:\cYY_{c} \to \cYY_{\tilde{c}}$ which
is equivariant with respect to $H$. By construction of $\sigma$ and
$\tau$, it suffices to find $i_0: \cZZ_{c} \to \cZZ_{\tilde{c}}$ which is
equivariant with respect to $\sigma$ and $\varphi_1$, and such that
the quotient isomorphism $\widehat{i}_0: \widehat{\cZZ}_{c} \to 
\widehat{\cZZ}_{\tilde{c}}$ is equivariant with respect to $\tau$. 
\par
Denote by $i: \cYY_{c} \to \cYY_{\tilde{c}}$ the canonical isomorphism and
 try $i_0 := \varphi^j \circ i$, for a suitably chosen $j$.
Then $i_0$ is automatically $\varphi_1$-equivariant.
Let $\pi_1$ (resp.\ $\pi_2$) denote the maps
from $C$ to the intermediate cover given by $s^n=t$ (resp.\
 $s^m = (t-1)$). By hypothesis we have $\pi_1(c) = \zeta_n^{e_1}
\pi_1(\tilde{c})$  and $\pi_2(c) = \zeta_m^{e_2} \pi_2(\tilde{c})$,
where $\zeta_m$ (resp.\ $\zeta_n$) is an $m$-th (resp.\ $n$-th)
root of unity.
We have 
$$ \tau \circ \widehat{i} = \varphi^{2me_1} \circ \widehat{i} \circ \tau, \quad
\sigma\circ i = \varphi^{2ne_2+2me_1} \circ i \circ \sigma. $$
The equivariance properties for $i_0 = \varphi^j
\circ i$ impose the condition
$$(\alpha -1)j + 2me_1 \equiv 0 \mod N/\gamma, \quad \text{and} \quad
-2j + 2ne_2 + 2me_1 \equiv 0 \mod N$$
on $j$. 
These conditions are equivalent to 
$$ -2j + 2me_1 \equiv 0 \mod 2n/\gamma \quad  \text{and} \quad
 -2j + 2ne_2 \equiv 0 \mod 2m.$$
We can solve  $j$, since $\gcd(m,n/\gamma)=1$.
\par
To see that the projective affine group is not larger than 
$\Delta(m,n,\infty)$ 
for $m\neq n$ we note that a larger projective affine group is again 
a triangle group. Singerman (\cite{Si72}) shows that any inclusion 
of triangle groups
is a composition of inclusions in a finite list. The 
case $\Delta(m,m,\infty) \subset \Delta(2,m,\infty)$ 
is the only one case that might occur here.
\hspace*{\fill} $\Box$
\par
\begin{Cor}
All $(m,n,\infty)$-triangle groups for $m,n >1$ and $mn \geq 6$ 
arise as projective affine groups of translation surfaces with
$\Delta(m,m,\infty)$ as possible exception.
\end{Cor}
\par
We determine the basic geometric invariant of the Teichm\"uller
curves constructed in Theorem \ref{hisTeich}.
\par
\begin{Thm} \label{signmninf}
In case S and DE the generating differential $\omega_X$ has $\gamma/2$
zeros and in the cases O and OE the generating differential $\omega_X$ has
$\gamma$ zeros.
\end{Thm}
\par
{\bf Proof:} We only treat the cases O and OE. The cases S and DE are similar.
We calculate the zeros of the pullback $\omega_Y$
of $\omega_X$ to the corresponding fiber $Y$ of $\cYY$. The differential 
$\omega_i$  on $Z$ has zeros of order ${a_1(i)\gamma}/{\gamma_1}-1$ 
(resp.\ 
${a_2(i)\gamma}/{\gamma_1}-1$) at the $\gamma_1/\gamma$
points over $0$ (resp.\ $1$). It has zeros of order
${a_3(i)\gamma}/{\gamma_2}-1$ (resp.\ ${a_4(i)\gamma}/{\gamma_2}-1$) 
at the $\gamma_2/\gamma$ points over $t$ (resp.\ $\infty$).
Therefore, the pullback of $\omega_i$ to $Y$ has zeros of order $a_\mu(i) -1$ 
at the $\gamma$ preimages of $\mu=0, 1,t,\infty$.
The differential $\omega_Y$ is a linear combination with non-zero
coefficients of $\omega_1$, $\omega_{-1}$ and two differentials that 
are pulled back from $Z^\tau$. The vanishing orders of these differentials on 
$Z^\tau$
are obtained from those of $\omega_1$ and $\omega_{-1}$ on $Z$ 
by replacing $a_1$ by $a_4$, $a_2$ by $a_3$,
and conversely. Since the $a_\mu$ are pairwise distinct, we
conclude that $\omega_Y$ vanishes at the (in total) $4\gamma$
preimages of $\{0,1,t,\infty\}$ of order
$\min\{a_1,a_2,a_3,a_4\} -1 = a_4-1$. Since $\omega_Y$ vanishes
also at the $4m+4n$ ramification points of $Y \to X$ we deduce
that it vanishes there to first order and nowhere else.
The $4\gamma$ zeros at the non-ramification points yield
the $\gamma$ zeros of $X$.
\hspace*{\fill} $\Box$
\par

\subsection{Comparison with Ward's results}\label{Wardsec}
In this section we compute an explicit equation for one particular
fiber of the family $\ol{f}=\ol{f}(m,n):\overline{\cXX}\to \overline{C}$. This
fiber, $\cXX_c$, is chosen such that $\cXX_c$ is a cyclic cover of a
projective line. This result is used in \S \ref{Billiards} to
realize $\cXX_c$ via unfolding a billiard table, for small $m$. In
this section we show moreover that for $m=3$, the family
$\ol{f}:\overline{\cXX}\to \overline{C}$ coincides with the family of
curves constructed by Ward \cite{Wa98}.

\par The assumptions on $m$ and $n$ in the following theorem are not
necessary. We include them to avoid case distinctions.  The reader can
easily work out the corresponding statement in the general situation,
as well.  We use the same notation as in the rest of this section. In
particular, $\ol{\pi}:\overline{C}\to \PP^1_t$ denotes the natural
projection of $\overline{C}$ to the $t$-line defined in Notation \ref{pinot}. 
One
may of course interchange the role of $m$ and $n$ in the theorem. In
that case one should consider the fiber of $\cXX$ in a point of $\overline{C}$
above $t=1$, instead.  \par
\begin{Thm} \label{Wardexplicit}
Suppose that $m$ and $n$ are relatively prime and $n$ is odd. Then a 
fiber of $\ol{f}:\ol{\cXX} \to \ol{C}$ over a point of $\ol{\pi}^{-1}(0)$ is a
$2n$-cyclic cover of the projective line branched at
$(m+3)/2$ points if $m$ is odd and $(m+4)/2$ points
if $m$ is even. 
\begin{itemize}
\item[(a)]
For $m$ odd this cover is given by the equation
$$ X_0: \qquad y^{2n} = (u-2) \prod_{k=1}^{(m-1)/2} \left(u-
2\cos\left(\frac{2k\pi}{m}\right)\right)^{2}.$$ 
The generating differential of the Teichm\"uller curve is
$$  \omega_0 = \frac{y \,\dd u}{(u-2)\prod_{k=1}^{(m-1)/2} (u- 
2\cos(2k\pi/m))}.$$
\item[(b)]
For $m$ even this cover is given by the equation
$$ X_0:\qquad y^{2n} = (u-2)^n \prod_{k=1}^{m/2} \left(u- 2\cos 
\left(\frac{(2k-1)\pi}{m}\right)\right)^{2}.$$
The generating differential of the Teichm\"uller curve
is 
$$ \omega_0 = \frac{y \,\dd u}{(u-2)\prod_{k=1}^{m/2} 
\left(u- 2\cos((2k-1)\pi/m)\right)}.$$
\item[(c)]
For $m=3$ the  surface $(X_0, \omega_0)$ is
the translation surface found by Ward.
\end{itemize}
\end{Thm}

{\bf Proof:} Our simplifying assumptions imply that $\gamma=1$ and
$\cZZ \cong \cYY$ (Remark \ref{simplerem}). Let $c$ be a point of
$\ol{C}$ with $\ol{\pi}(c)=0$. Then the fiber $\cYY_c$ of $\cYY$
consists of two isomorphic irreducible components, $Y_0^j$, given by
the affine equation $y^N =x^{a_1}(x -1)^{a_3}$. Note that $Y^1_0\to
\PP^1_{x}$ is branched at $x=\infty$ of order $m$. The fiber
$X_0:=\cXX_c$ of $\cXX$ is the quotient of $ Y_0^1$ by $\rho$.

From the presentation of $G$ (Lemma \ref{autlem2}) we deduce that
$\varphi^k_1$ commutes with $\rho$ if and only if $k$ is a multiple of
$m$. We denote by $A$ the abelian subgroup of $\Aut(Y_0^1)$ generated
by $\rho$ and $\varphi_1^m$.  One computes that the quotient of
$Y_0^1$ by $\langle \varphi_1^m\rangle$ has genus zero. We denote this
quotient by $\PP^1_z$. Here $z$ is a parameter on $\PP^1_z$ such that
$\PP^1_z\to \PP^1_x$ is given by
\[
z^m=\left(\frac{x-1}{x}\right)^n.
\]
 Let $\PP^1_u$ be the quotient of $Y^1_0$ by
$A$. The subscript $u$  denotes a coordinate which is defined
below. We obtain the following diagram of covers
\[
\xymatrix{ & Y_0^1 \ar[dr]^{\langle \varphi_1^m \rangle}_{p_2}
\ar[dl]_{\langle \rho \rangle}^{p_1}\\
 X_0 \ar[dr]_{q_1} &&\PP^1_z
\ar[dl]^{q_2} \ar[dr]^{\langle \varphi_1 \bmod m \rangle}\\ 
& \PP^1_u & & \PP^1_{x}  }
\]

Suppose that $m$ is odd. After replacing $y$ by $z^{(n+1)/2}/y$, we
find that $Y^1_0\to \PP^1_z$ is given by
\begin{equation}\label{Y10eq}
y^{2n}=\frac{(z^m-1)^2}{z^m}.
\end{equation}
Here we use that $n$ is odd.
Recall that $\rho\in \Aut(\PP^1_x)$ is given by $\rho(x)=1-x$.  It
follows from Lemma \ref{autlem1} that $\rho$ lifts to an automorphism
of order $2$ of $\PP^1_z$ which has one fixed point in the fiber above
$x=1/2$. Without loss of generality, we may assume that $\rho(z)=1/z$.
Therefore, $u:=z+1/z$ is an invariant of $\rho$; it is a
parameter on $\PP^1_u$. We find an equation for $X_0$ by rewriting
(\ref{Y10eq}) in terms of $y$ and $u$. Noting that
\[
u-(\zeta_m^i+\zeta_m^{-i})=\frac{(z-\zeta_m^i)(z-\zeta_m^{-i})}{z},
\]
we find the equation in (a).  The differential form $\omega_0$ in (a)
 is a holomorphic differential form with a zero only in
 $u=0$. Therefore Theorem \ref{signmninf} implies that $\omega_0$ is a
 generating differential form. 

Specializing to $m=3$, we find the
 equation found by Ward (\cite{Wa98} \S5). This proves (c).

Suppose now that $m$ is even.  After replacing $y$ by $z^{m/2}/y$, we
find that $Y^1_0\to \PP^1_x$ is given by
\[
y^{2n}=\frac{(z^m-1)^2}{z^{m+n}}.
\]
In this case the automorphism $\rho$ of $\PP^1_x$ lifts to an
automorphism of $\PP^1_z$ with two fixed points in the fiber above
$x=1/2$. Without loss of generality, we may suppose that
$\rho(z)=\zeta_m/z$.  Therefore, $u:=\zeta_{2m}^{-1}z+\zeta_{2m}/z$ is
an invariant of $\rho$ which we regard as parameter on $\PP^1_u$.
Here $\zeta_{2m}$ is a square root of $\zeta_m$. One computes that
\[
u-(\zeta_{2m}^{2i-1}+\zeta_{2m}^{-2i+1})=\zeta_{2m}^{-1}
\frac{(z-\zeta_m^i)(z-\zeta_m^{1-i})}{z}, \quad
u-2=\zeta_{2m}^{-1}\frac{(z-\zeta_{2m})^2}{z}.
\]
After replacing $y$ by $c(z-\zeta_{2m})/y$ for a suitable root of
unity $c$, we find the equation for $X_0$ which is stated in (b). The
expression for $\omega_0$ follows as in the proof of (a).
\hspace*{\fill}$\Box$

\par
\section{Primitivity} \label{primsection}
A {\em translation covering} is a covering $q: X \to Y$ between translation
surfaces $(X,\omega_X)$ and $(Y,\omega_Y)$ such that $\omega_X = q^*
\omega_Y$. A translation surface $(X,\omega_X)$ is called {\em geometrically
  primitive} if it does not admit a translation covering to a surface $Y$ with
$g(Y) < g(X)$.  \newline A Veech surface $(X,\omega)$ is called {\em
  algebraically primitive} if the degree $r$ of the trace field extension over
$\QQ$ equals $g(X)$. Algebraic primitivety implies geometric primitivety, but
the converse does not hold (\cite{Mo04b}). In loc.\ cit.\ Theorem\ 2.6 it is
shown that a translation surface of genus greater than one covers a unique
primitive translation surface.
\par
Obviously the Veech examples ($p: \cUU \to \widetilde{C}$ in the notation of
Theorem \ref{Teichnevenquottau}) for $n=2\ell$ and $\ell$ prime and those for
$(2,n,\infty)$ (compare to Remark \ref{oddVeech}) are algebraically primitive.
We will not give a complete case by case discussion of primitivity of the
$(m,n,\infty)$-Teichm\"uller curves, but restrict to the case that $m$ and $n$
are odd and relatively prime.  Comparing the degree of the trace field
$[\QQ(\zeta_m+\zeta_m^{-1}, \zeta_n+\zeta_n^{-1}:\QQ] = r \leq
\phi(m)\phi(n)/4$ with the genus  (Corollary \ref{genuscor}), we deduce
that the fibers of $\cXX \to C$ are never algebraically primitive. 
  Nevertheless, we  show that
 there are infinitely many geometrically primitive
ones.
\par
\begin{Thm} \label{mnisprim}
Let $m$, $n$ distinct odd primes. Then the Veech surfaces arising 
from the $(m,n,\infty)$-Teichm\"uller curve $f:\cXX \to C$ of
Theorem\ \ref{hisTeich} are geometrically primitive.
\end{Thm}
\par
{\bf Proof:} Let $(X, \omega_X)$ be such a Veech surface and suppose there is
a translation covering $q: X \to Y$.  Then $g(Y) \geq r$, by \cite{Mo04b}
Theorem\ 2.6.  Theorem \ref{signmninf} implies that the generating
differential has only one zero on $\cXX_c$. Therefore the cover $q$ is totally
ramified at this zero, and nowhere else. This contradicts the Riemann--Hurwitz
formula.  Namely, a degree two cover cannot be branched in exactly one point.
If the degree $d$ of $q$ is larger than $2$, we obtain a contradiction with
$g(Y) \geq r$. \hspace*{\fill} $\Box$
\par
\begin{Rem}\label{series}
{\rm
At the time of writing the authors are aware of the
following series of examples of Teich\-m\"uller curves: 
The triangle constructions in \cite{Ve89} and \cite{Wa98}
and the Weierstrass eigenform or Prym eigenform constructions
in \cite{McM03} and \cite{McM05}. Besides them there
is a finite number of sporadic examples.
}\end{Rem}
\par
\begin{Cor} 
The Veech surfaces arising from the case $(m,n,\infty)$
with $m$, $n$ sufficiently large distinct primes
are not translation covered by any of the Veech surfaces 
listed is Remark \ref{series}.
\end{Cor}
\par
{\bf Proof:} Recall that translation coverings between
Veech surfaces preserve the affine group up to commensurability.
In particular, they preserve the trace field.
\newline
Choose $m$ and $n$ sufficiently large such that
the trace field $K$ of the $(m,n,\infty)$-triangle group
is none of the trace fields occurring
in the sporadic examples and such that the genus is
of the Veech surface is larger than $5$. This implies that the
surface cannot be one of examples in \cite{McM03} and \cite{McM05}.
There is only a finite list of arithmetic triangle groups (\cite{Ta77}).
We choose $m>3$ and $n>5$ such that $K$ is not one of the
trace fields in this finite list. Non-arithmetic
lattices have a unique maximal element (\cite{Ma91})
in its commensurability class and the $(m,n,\infty)$-triangle
groups are the maximal elements in their classes. Since the
$(2,n,\infty)$- and $(3,n,\infty)$-triangle groups are the
maximal elements in the commensurability classes of the
examples of \cite{Ve89} and \cite{Wa98}, these examples 
cannot be a translation cover of the examples given
by Theorem\ \ref{hisTeich} for $(m,n)$ chosen as above.
\hspace*{\fill} $\Box$
\par
\begin{Rem} {\rm
Even in the cases that the Veech surfaces with affine group
$\Delta(m,n,\infty)$ are geometrically primitive, 
Theorem 2.6 of \cite{Mo04b}
does not exclude that there are other primitive Veech surfaces
with the same affine group. Such examples are provided by
Theorem 3' of \cite{HuSc01} for $n=\infty$. By Remark \ref{VHSunique}
we know a rank $2r$ subvariation of Hodge structures of the
family of curves generated by such a Veech surface. In particular,
we know  $r$ of the Lyapunov exponents $\lambda(\LL_i)$.
}\end{Rem}

\section{Billiards} \label{Billiards}
In this section we approach Teichm\"uller curves uniformized by
triangle groups in the way Veech and Ward did in \cite{Ve89} and
\cite{Wa98}. We start by presenting two series of billiard tables
$T{(m,n,\infty)}$, for $m=4,5$. These tables are (rational) $4$-gons
in the complex plane.  We show that the affine group of the
translation surface $X(m,n,\infty)$ attached to $T{(m,n,\infty)}$ is
the $(m,n,\infty)$-triangle group, for $m=4,5$. This part is
independent of the previous sections, and requires only elementary
notions of translation surfaces (see \cite{MaTa02} or \S
\ref{Teich}).  The proof we give that these billiard tables define
Teichm\"uller curves is combinatorically more complicated than the
analoguous proof for the series of Teichm\"uller curves found by Veech
and Ward. This suggests that it would have been difficult to find
these billiards by a systematic search among $4$-gons.

In \S \ref{coincide} we relate these explicitly constructed
billiard tables to our main realization result (Theorem
\ref{hisTeich}). Denote by $f=f(m,n):\cXX\to C$ the family of curves
constructed in \S \ref{realize}. This family defines a finite map
from $C$ to $M_g$, for a suitable integer $g\geq 2$. The image of this
map is a Teichm\"uller curve whose (projective) affine group is the
$(m,n,\infty)$-triangle group. We have shown in Theorem
\ref{Wardexplicit} that a suitable fiber $\cXX_c$ of $\cXX$ is a
$2n$-cyclic cover of the projective line which we described
explicitly. In this situation, one may use a result of Ward to find
the corresponding billiard table $T[m,n,\infty]$. We show that
$T[m,n,\infty]$ may be embedded in the complex plane (i.e.\ without
self-crossings) if and only if $m\leq 5$. For $m=2, 3$ we find back
the billiard tables found by Veech and Ward. We show that the tables
we obtain for $m=4,5$ are the ones we already constructed.

Consider a compact polygon $P \subset \RR^2 \cong \CC$ in the plane
whose interior angles are rational multiples of $\pi$. The linear
parts of reflections along the sides of the polygon generate a finite
subgroup $G \subset O_2(\RR)$.  If $s$ is a side
of $P$ we write $\sigma_s$ for the linear part of the reflection in
the side $s$. One checks that for sides $s$ and $t$ of $P$ we have
\[
\sigma_{\sigma_s(t)}=\sigma_s\sigma_t\sigma_s.
\]

We define an equivalence relation on $G$ as follows. We write
$\sigma_1\sim \sigma_2$ if the reflected polygons $\sigma_1(P)$ and
$\sigma_2(P)$ differ by a translation in $\CC$. 
Let $G_0\subset G$ represent the 
equivalence classes of this relation.  By gluing copies of $P$ we
obtain a compact Riemann surface
$$ X = \left(\coprod_{g \in G_0} gP \right)/ \approx, $$ where $\approx$
denotes the following identification of edges: if $gP$ is obtained
from $\tilde{g}P$ by a reflection $\sigma$ along a side $s$ of
$\tilde{g}P$, then $s$ is glued to the edge $\sigma(s)$ of $gP$ by a
translation.  \par The holomorphic one-from $\dd z$ on $P$ and its
copies is translation invariant, hence defines a $1$-form $\omega$ on
$X$.  We say that the translation surface $(X,\omega)$ is {\em
obtained by unfolding $P$}. The trajectories of a billiard ball on $P$ 
correspond to straight lines on $X$.  In \cite{McM03} $X$ is called
the {\em small surface} attached to $P$.  The translation surface has
a finite number of points where the total angle exceeds $2\pi$. These
are called {\em singular points}. They correspond to the
zeros of $\omega$.  \par

\subsection{The tables $T(5,n,\infty)$}\label{5nsec}
\begin{figure}[h] 
\centering
\includegraphics[scale=0.65]{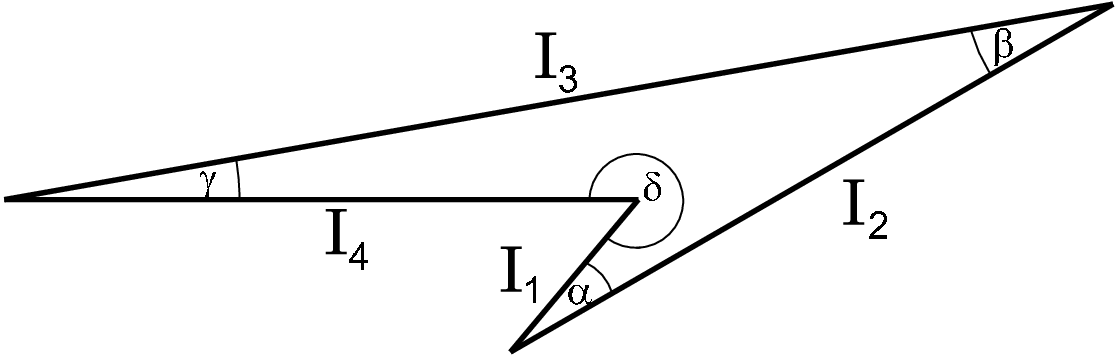}
\caption{Billiard table $T{(5,n,\infty)}$, for $n=9$}\label{billiard59}
\end{figure}

Let $n \geq 7$ be an odd integer which is not divisible by $5$.  We
define a billiard table $T{(5,n,\infty)}$ as follows (Figure
\ref{billiard59}). The billiard table $T(5,n,\infty)$ is a $4$-gon in
the complex plane with angles $\alpha=\beta=\pi/n$ and
$\gamma=\pi/2n$, as indicated in the picture. 
We denote by $I_1, \ldots, I_4$ the vectors corresponding to the sides
of the polygon which we regard as complex numbers. We rotate and scale
the billiard table such that $I_4 = 1$ and
$$|\Re(I_3)| = \cos(\pi/n) + \cos(\pi/5).$$ In particular, $I_4$
points in the direction of the positive $x$-axis. This determines the
table uniquely.  \par
\begin{figure}[h]
\centering
\includegraphics[scale=0.79]{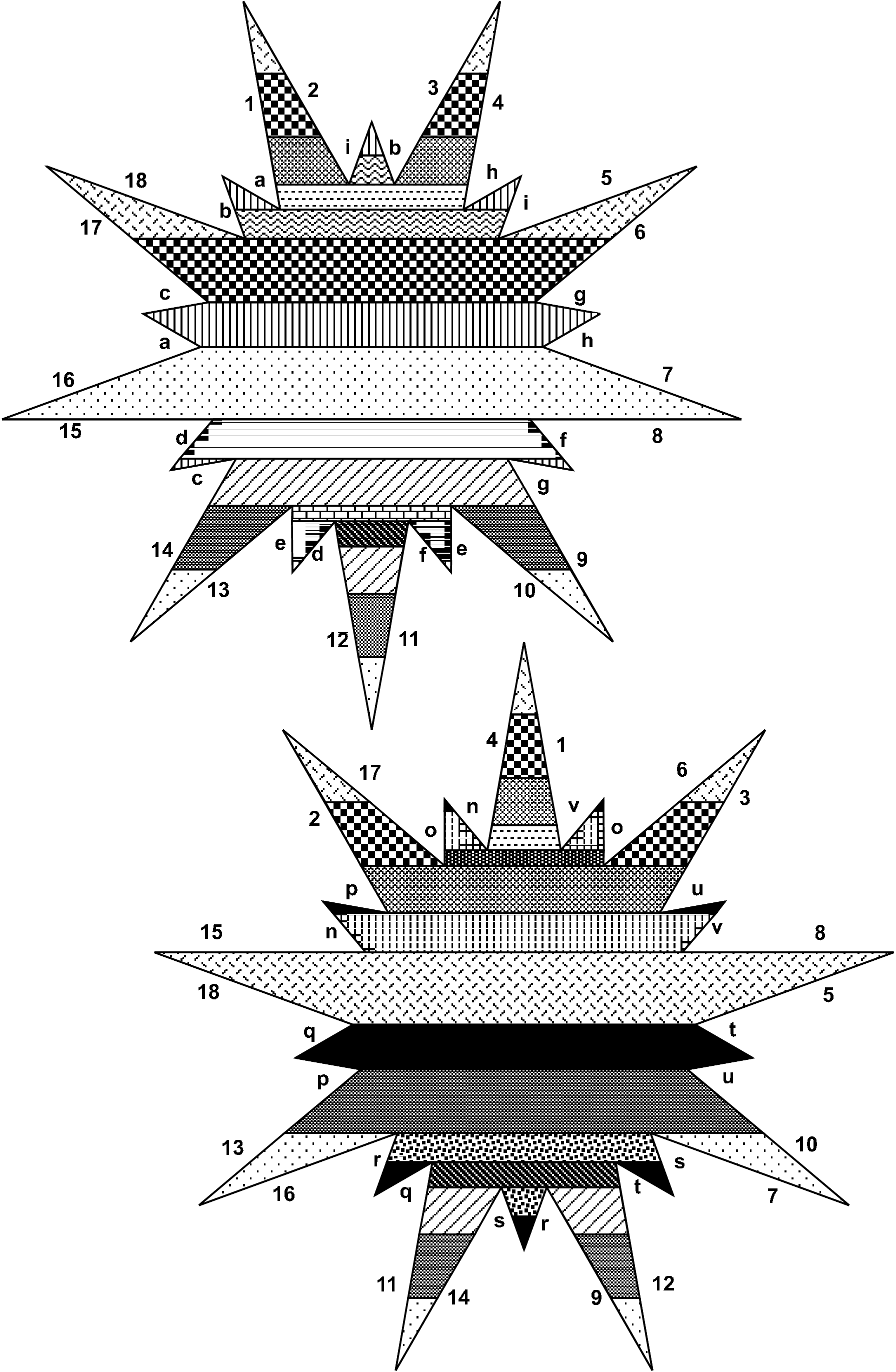}
\caption{Cylinder decomposition of $X(5,9,\infty)$}  \label{cyldec}
\end{figure}

 We now construct the translation surface
obtained by unfolding the table $T{(5,n,\infty)}$ (Figure
\ref{cyldec}).  Reflecting the billiard table $2n$ times in the
(images of the) sides $I_2$ and $I_3$ yields the upper star; it
consists of alternating long and short points. The second star is
obtained by reflecting the first star in the side $I_4$ of the
billiard table (this is the side marked by $15$ in Figure
\ref{cyldec}).  The two stars can be glued together to a translation
surface $X:=X(5,n,\infty)$: sides denoted by the same  letters
or numbers are glued by translations. Note that the tips of the
`long points' (resp.\ the `short points') of the stars correspond
to one point of the translation surface $X$; both points of $X$ are
not singularities, since the total angle is $2\pi$. There is one
singularity; it corresponds to the angle $\delta$.  The genus of $X$
is $g=2(n-1)$.
\par
\begin{Thm} \label{59genTrGrp}
Let $n \geq 7$ be odd and not divisible by $5$. Then the affine group of 
$X(5,n,\infty)$ contains the elements 
$$R = \left(\begin{array}{cc} \cos(\pi/n) & -\sin(\pi/n) \\
\sin(\pi/n) & \cos(\pi/n) \end{array}\right) \quad \text{and} \quad
T= \left(\begin{array}{cc} 1 & 2\frac{\cos(\pi/n)+\cos(\pi/5)}{\sin(\pi/n)} \\
0& 1 \end{array}\right).$$
The elements $R,T  \in \PSL_2(\RR)$ generate the Fuchsian
triangle group $\Delta(5,n,\infty)$. In particular,
$X(5,9,\infty)$ is a Veech surface.
 \end{Thm}
\par
{\bf Proof:} Rotation around the center of the stars defines an
affine diffeomorphism of the surface $X(5,n,\infty)$.  Its  derivative  is $R$.
\par
We rotate $X(5,n,\infty)$ as in Figure \ref{billiard59} and 
Figure \ref{cyldec}, i.e.\ such that the edge $I_4$ resp.\ the one
with label $15$ is horizontal and to the left of the center of the star.
\clearpage

We now consider the horizontal foliation defined by $\omega$. Recall
that a {\em saddle connection} is a leaf of the foliation that begins
and ends in a singularity.  In a dense set of directions, the saddle
connections divide $X$ into {\em metric cylinders}, see for example
\cite{MaTa02}, \S 4.1.  We claim that in the horizontal direction
$X$ decomposes into $g=2(n-1)$ metric cylinders. We distinguish two
types of cylinders. Each cylinder corresponds to one  shading style in 
Figure~\ref{cyldec}. 

The cylinders of type $1$, denoted by $C_i$, are those that are glued
together from pieces from both stars. An example is the checkered
cylinder. Since the second star is obtained from the first by
reflection, the cylinders $C_i$ appear in pairs, as can be seen from
Figure \ref{cyldec}. There is a bijection between the cylinders of
type $1$ and pairs of long points. For example, 
the checkered cylinder
corresponds to the long points $17$-$18$ and $5$-$6$. Here a `pair'
consists of an orbit of length $2$ of long points under the reflection in the
vertical axis. The two vertical long points correspond to orbits of
length one, and hence do not correspond to a cylinder of type $1$. We
conclude that the number of cylinders of type $1$ is $n-1$.

The cylinders of the type $2$, denoted by $\widetilde{C}_i$, are those
that consist of pieces of one star. An example is the black
cylinder. These cylinders also come in pairs. There is a bijection
between cylinders of type $2$ and pairs of short points. Therefore the
number of cylinders of type $2$ is also $n-1$.

   The width and the  height  of a pair of cylinders of type $1$,  for an
appropriate numbering, is given  by
\begin{equation}\label{wh}
\begin{split}
 w_k &= 2 |I_3| \cos \frac{(n-2k)\pi}{2n} \quad \text{and} \\ h_k &=
|I_4|\left(\sin\frac{(n+1-2k)\pi}{ 2n} - \sin \frac{(n-1-2k)\pi}{2n}
\right) = 2|I_4| \sin \frac{\pi}{2n} \cos\frac{(n-2k)\pi}{2n}
\end{split}
\end{equation}
for $k=1,\ldots,(n-1)/2$. This is seen by cutting the points of the
stars into pieces, and translating these pieces so that one obtains
$2(n-1)$ connected cylinders, one for each shading style. One then uses the
rotation and reflection symmetries of the original star.  \par
Similarly, the widths and heights of pairs of cylinders
$\widetilde{C}_i$, for an appropriate numbering, are given by
\begin{equation}\label{wht}
\begin{split}
 \widetilde{w}_k &= 2|I_2| \cos\frac{(n-2k)\pi}{2n} \quad
\text{and} \\
 \widetilde{h}_k &=  
|I_1|\left(\sin\frac{(n+2-2k)\pi}{
2n} - \sin \frac{(n-2-2k)\pi}{2n} \right) =
2|I_1|
\sin \frac{2\pi}{2n} \cos\frac{(n-2k)\pi}{2n}
\end{split}
\end{equation}
for $k=1,\ldots,(n-1)/2$.
\par
The moduli of these cylinders are 
$$m_k := h_k/w_k = |I_4| \sin \frac{\pi}{2n}/ |I_3| 
\quad \text{and} \quad\widetilde{m}_k := \widetilde{h}_k/
\widetilde{w}_k = |I_1| \sin \frac{2\pi}{2n}/ |I_2|. 
$$ Note that $m_k$ and $\widetilde{m}_k$ are independent of $k$.  \par
We claim that $m_k / \widetilde{m}_k = |I_4||I_2|\sin(\pi/n)/
|I_3||I_1| \sin(\pi/2n)=1$, that is that the moduli of all the
cylinders are identical. This is equivalent to
\begin{equation}\label{modulieq}
\frac{|I_2|}{|I_1|}=\frac{|I_3|}{|I_4|}\frac{\sin(\pi/n)}{\sin(\pi/2n)}.
\end{equation}
Since we assumed that $I_4=1$, the right hand side is equal to $2|\Re(I_3)|$.

Using  the geometry of the billiard table one shows that
\[
\begin{split}
      |I_2|\cos(3\pi/2n) & -  |I_1|\cos(5\pi/2n) 
  =  |\Re(I_3)|-|I_4|=|\Re(I_3)|-1,\\ 
   |I_2|\sin(3\pi/2n) & - 
  |I_1|\sin(5\pi/2n)  =  \Im(I_3) = \Re(I_3) \tan(\pi/2n).  
\end{split}
\]
This implies that
\begin{equation}\label{tableeq}
\frac{|I_2|}{|I_1|} = \frac{-(\Re(I_3)-1)\sin(5\pi/n) + \Re(I_3) 
\tan(\pi/2n) \cos(5\pi/n)}{-(\Re(I_3)-1)\sin(3\pi/n) + \Re(I_3) 
\tan(\pi/2n) \cos(3\pi/n)}.
\end{equation}

The minimal polynomial of $\Re(I_3)$ over $\QQ(\cos(\pi/n))$ is
\begin{equation}\label{minpolyeq}
 X^2 - (2\cos(\pi/n) +1/2) X + (\cos^2(\pi/n) + \cos(\pi/n)/2 -
1/4).
\end{equation}
 One deduces (\ref{modulieq}) from (\ref{tableeq}), (\ref{minpolyeq})
and the addition laws for sines and cosines.

From the claim (\ref{modulieq}), we deduce that $T$ is contained in the affine
group of $X$. Namely, fixing  the horizontal lines and postcomposing
local charts of the cylinders by $T$ defines an affine diffeomorphism
whose derivative is $T$ (compare to \cite{Ve89} Proposition 2.4 or
\cite{McM03} Lemma 9.7).  \par It remains to prove that $ R$ and $T$
generate the $(5,n,\infty)$-triangle group.  One constructs the
hyperbolic triangle in the extended upper half plane with corners
$i\infty, i$ and
$$z_0 =  \frac{\cos (\pi/n) + \cos(\pi/5)}{\sin(\pi/n)} + i \,
\frac{\sin(\pi/5)}{\sin(\pi/n)}$$
bounded by the vertical
axes through $i$ and $z_0$ 
and the circle around $\cot(\pi/n)$ with radius $1/\sin(\pi/n)$.
The interior angles at $i$ and $z_0$ are indeed $\pi/n$ and $\pi/5$.
By Poincar\'e's Theorem this triangle plus its reflection along the 
imaginary axis is a fundamental domain for the group generated by
$R$ and $T$. 
\par
The last claim follows now from the standard criterion
to detect Teichm\"uller curves, see e.g.\ \cite{McM03} Corollary 3.3. 
 \hspace*{\fill} $\Box$
\par
\begin{Rem} \label{numofcyl}
  {\rm Assuming the comparison results which will be proved in Section
    \ref{coincide} below, the number of cylinders in, say, the horizontal
    direction is already determined by results of the previous sections.
\par
Consider the family of translation surfaces $\diag(e^t, e^{-t})\cdot 
(X_0, \omega_0)$, where $(X_0,\omega_0)$ is as in Theorem \ref{Wardexplicit}.
This family converges for $t \to \infty$ to a singular fiber of 
$\ol{f}: \ol{\cXX} \to \ol{C}$ and by \cite{Mr75}  the number of 
cylinders in the  horizontal direction equals the number of singularities 
of the singular fiber $X_\infty$.
\par
Since all the local systems $\LL_i$ as in the proof of
Theorem \ref{hisTeich} have non-trivial parabolic monodromy
around points in $\ol{\pi}^{-1}(\infty)$ the arithmetic genus of $X_\infty$
is zero. Since $\omega_0$ has only one zero, $X_\infty$ is irreducible 
and hence the number of singularities of the fiber
$X_\infty$ equals $g(X_\infty)$. 
}\end{Rem}
\par
\subsection{The tables $T(4,n,\infty)$}\label{4nsec}
\begin{figure}[h] 
\centering
\includegraphics[scale=0.65]{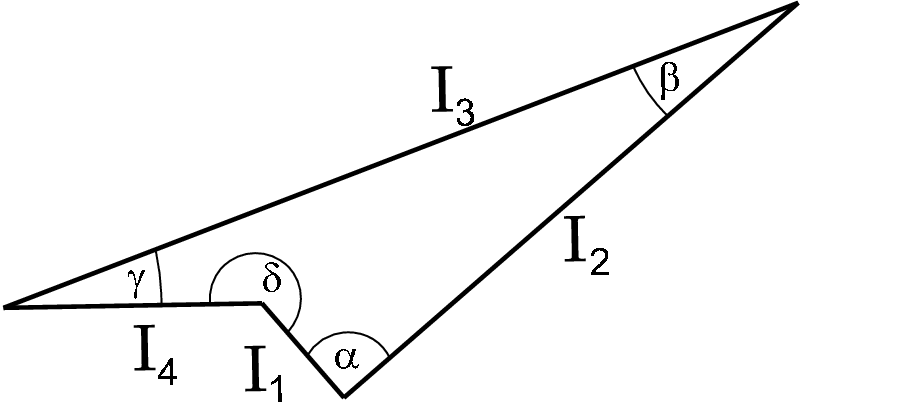}
\caption{Billiard table $T{(4,n,\infty)}$, for $n=9$}\label{billiard49}
\end{figure}

Let $n \geq 5$ be odd. We define a billiard table $T{(4,n,\infty)}$ as
follows. The billiard table is a again a $4$-gon in the complex plane
with angles $\alpha=\pi/2$ and $\beta=\gamma=\pi/n$, as indicated in
Figure \ref{billiard49}. We denote by $I_1, \ldots, I_4$ the vectors
corresponding to the sides of the polygon. We regards these vectors as
complex numbers. 
We scale and rotate the billiard table such that $I_4 = 1$ and such
 that $|I_3| = 2(\cos(\pi/n)+\cos(\pi/4)).$ This determines the table
 uniquely.  

\begin{figure}[h] 
\centering
\includegraphics[scale=0.45]{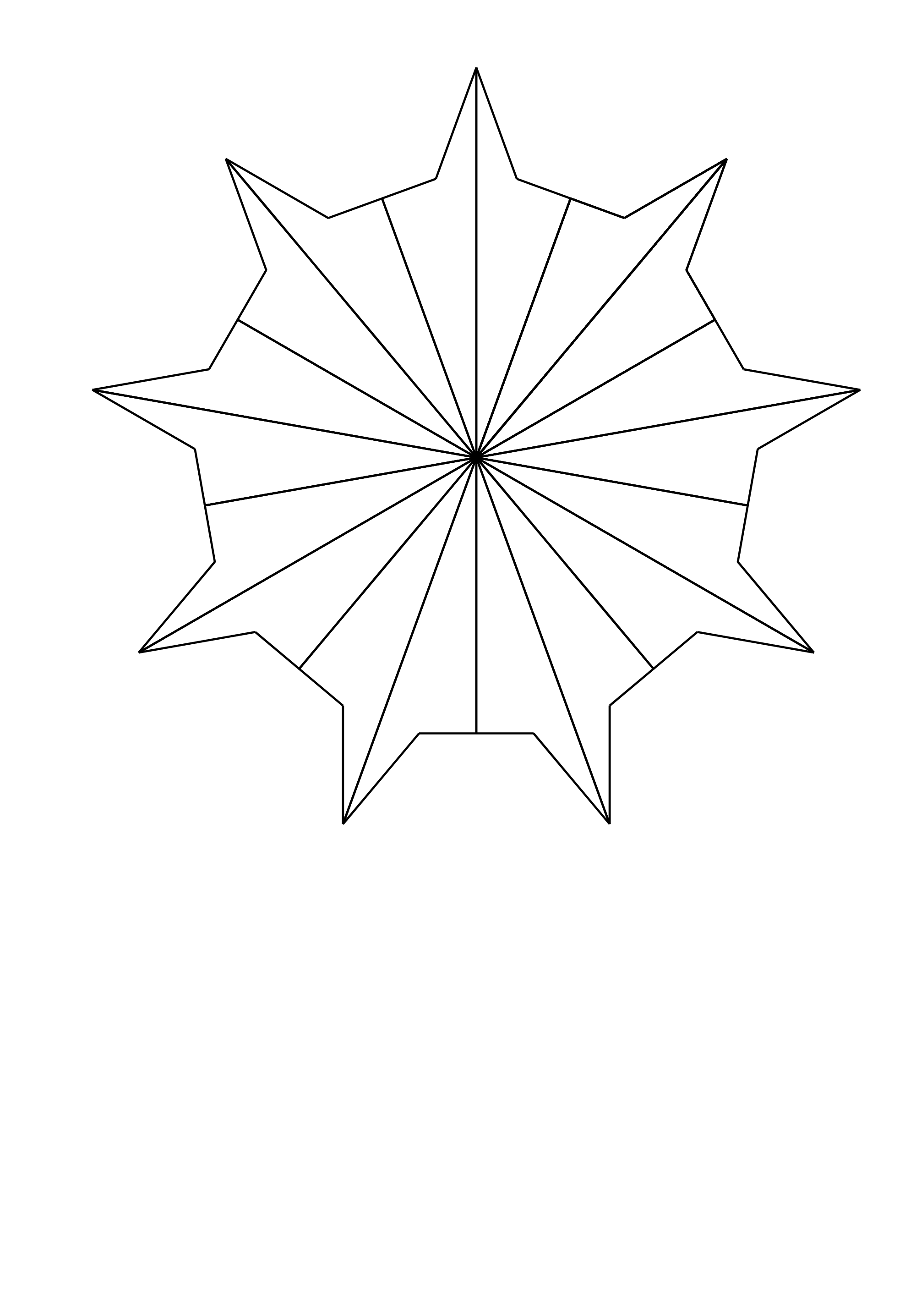}
\caption{Billiard table $T{(4,9,\infty)}$, partially unfolded}
\label{49unfolded}
\end{figure}

 The translation surface $X:=X(4,n,\infty)$ obtained by unfolding
 $T{(4,n,\infty)}$ looks similar to the one obtained from
 $T{(5,n,\infty)}$. It can be obtained by identifying parallel sides
 of two stars. The first star is illustrated in Figure
 \ref{49unfolded}. The second star is obtained from the first by
 reflection in the horizontal axis. 
The translation surface $X(4,n,\infty)$ has one singularity which
corresponds to the vertex of the billiard table with angle $\delta$.
Its genus is $g=3(n-1)/2$.  \par
\begin{Thm} \label{49genTrGrp}
Let $n\geq 5$ be odd. Then the affine group of 
$X(4,n,\infty)$ contains the elements 
$$R = \left(\begin{array}{cc} \cos(\pi/n) & -\sin(\pi/n) \\
\sin(\pi/n) & \cos(\pi/n) \end{array}\right) \quad \text{and} \quad
T= \left(\begin{array}{cc} 1 & 2\frac{\cos(\pi/n)+\cos(\pi/4)}{\sin(\pi/n)} \\
0& 1 \end{array}\right).$$
The elements $R,T  \in \PSL_2(\RR)$ generate the Fuchsian
triangle group $\Delta(4,n,\infty)$. In particular,
$X(4,9,\infty)$ is a Veech surface.
 \end{Thm}
\par {\bf Proof:} Rotation around the center of each of the stars
defines an affine diffeomorphism of $X(4, n, \infty)$ whose
derivative is $R$, as in the case $(5,n,\infty)$.

 We describe the cylinders in the horizontal direction. As for
$X(5,n,\infty)$, we distinguish two types of cylinders. The cylinders
of type $1$, denoted by $C_i$, are those that are glued together from
pieces of both stars. They correspond to pairs of sides which connect
two points. Here a pair of sides consist of two distinct sides which
are interchanged by reflection in the vertical axis. There are
$(n-1)/2$ such cylinders.  The widths and heights of these cylinders,
in an appropriate numbering, are given by
\[
 w_k = 4 |I_2| \cos\left(\frac{(n -2k)\pi}{2n}\right) \quad \text{and}
\quad h_k = 2 |I_1| \sin\left(\frac{k\pi}{n}\right), \qquad k=1,\ldots,(n-1)/2.
\]
 There are two cylinders with the same width and height, due to the symmetry.

The cylinders of type $2$, denoted by $\widetilde{C}_i$, are those
that consist of pieces of one star only. They correspond to pairs of
points of the stars. Here we use the same convention for pairs as
above. The number of such cylinders is also $(n-1)/2$. The widths and
heights of these cylinders are are
\[
\begin{split}
 \widetilde{w}_k &= 2 |I_3| \cos\left(\frac{(n -2k+2)\pi}{2n}\right)
 \quad \text{and} \\ \widetilde{h}_k &= 2 |I_4|
 \cos\left(\frac{(n -2k+2)\pi}{2n}\right)
 \sin\left(\frac{\pi}{n}\right), \qquad k=1,\ldots,(n-1)/2.
\end{split}
\]
The moduli of the cylinders are
$$ m_k := h_k/w_k=  |I_1|/2|I_2| \qquad \text{and}  \qquad \widetilde{m}_k=
\widetilde{h}_k/ \widetilde{w}_k = |I_4|\sin \frac{\pi}{n}/ |I_3|. $$
\par
As in the proof of Theorem~\ref{59genTrGrp},  one checks that
$m_k/\widetilde{m}_k= |I_1||I_3|/2|I_2||I_4|\sin\left( 
{\pi}/{n}\right) =1$, by using the geometry of the billiard table and
the minimal polynomial of $2(\cos(\pi/n)+\cos(\pi/4))$
over $\QQ(\cos(\pi/n))$.
The rest of the proof is analogous to the proof of Theorem~\ref{59genTrGrp}.
\hspace*{\fill} $\Box$
\par
\subsection{Comparison with Theorems \ref{hisTeich} and  \ref{Wardexplicit}} 
\label{coincide}
In this section we relate the billiard tables constructed in \S\S
\ref{5nsec} and \ref{4nsec} to the families of curves constructed in
\S \ref{realize}. For simplicity we suppose that $1<m<n$ are
relatively prime integers such that $n$ is odd. This assumption avoids
a case distinction. It is easy to work out the general statement. 

In Theorem \ref{hisTeich} we constructed a Teichm\"uller curve with projective
affine group $\Delta(m,n,\infty)$. We constructed a concrete finite cover,
$C$, of this Teichm\"uller curve. We denote by $\ol{C}$ the corresponding
projective curve. Over $\ol{C}$ there exists a universal family
$f:\ol{\cXX}\to \ol{C}$ of semistable curves.  In Theorem \ref{Wardexplicit}
we showed that there exist points $c$ of $\ol{C}$ such that the fiber
$X_0:=\cXX_c$ is a smooth curve which is a $2n$-cyclic cover of the projective
line branched at $(m+3)/2$ (resp.\ $(m+4)/2$) points if $m$ is odd (resp.\ 
even).  There also exist fibers of $\cXX$ which are $2m$-cyclic covers of the
projective line branched at $(n+3)/2$, but we do not regard these here since
it is convenient to have as few branch points as possible, for our purposes.
One may check that this is the most efficient way to represent a fiber of
$\cXX$ as an abelian cover of the projective line. This representation allows
us to use Schwarz--Christoffel maps (\cite{Wa98} Theorem C') to represent the
fiber $X_0$ of $\cXX$ as the unfolding of a billiard table, under certain
conditions (see below).

We first suppose that $m$ is odd. The $2n$-cyclic cover $X_0\to
\PP^1_u$ of Theorem \ref{Wardexplicit}.(a) is branched at the real
points $2\cos(2k\pi/m)$.
The Schwarz--Christoffel map is defined as
\[
SC(w) = \int_{0}^w (u-2)^{\frac{1}{2n}-1} \prod_{k=1}^{(m-1)/2} (u-
2\cos(2k\pi/m))^{\frac{1}{n} -1} \,\dd u.
\]
The integrand is the generating differential form $\omega_0$.  

The Schwarz--Christoffel map maps the real axis to a $(m+3)/2$-gon which we
denote by $T[m,n,\infty]$. If $T[m,n,\infty]$ has no self-crossings then $SC$
maps the upper half-plane bijectively to the interior of this $(m+3)/2$-gon.
The interior angles of $T[m,n,\infty]$ are $(m-1)/2$ times $\pi/n$ and once
$\pi/2n$, in this order.  The remaining angle is $2\pi-m\pi/2n \bmod 2\pi$
(resp.\ $\pi-m\pi/2n$) if $m\equiv 1\bmod{4}$ (resp.\ $m\equiv 3\bmod{4}$).
The number of self-crossings is therefore $(m-5)/4$ if $m\equiv 1\bmod{4}$ and
$(m-3)/4$ if $m\equiv 3\bmod{4}$. In particular, this number is zero if and
only if $m=3,5$.  For $m\geq 7$ it therefore unclear whether one can obtain
$(X_0, \omega_0)$ by unfolding a billiard table. However, it follows from our
results that one cannot do this via the usual theory of Schwarz--Christoffel
maps. Namely, for $m\geq 7$ one cannot represent a smooth fiber of $\ol{\cXX}$
as a cyclic cover of the projective line, such that the corresponding polygon
does not have self-crossings.

If $m=3$ or $5$, Theorem C' of \cite{Wa98} implies that the Veech surface
$(X_0, \omega_0)$ is obtained by unfolding the billiard table
$T[m,n,\infty]$. For $m=3$, we obtain Ward's family (compare to
Theorem \ref{Wardexplicit}.(c)). For $m=5$, the angles of the $4$-gon
$T[5, n, \infty]$ coincide with those of the billiard table
$T(5,n,\infty)$ which we constructed in \S \ref{5nsec}. We show
below that both $4$-gons are similar.

The case that $m$ even is analogous. The
Schwarz--Christoffel map
$$SC(w) = \int_{0}^w (u-2)^{\frac{1}{2}-1} \prod_{k=1}^{m/2} (u-
2\cos((2k-1)\pi/2m))^{\frac{1}{n} -1}\, \dd u$$ maps the real axis to
a $(m+4)/2$-gon $T[m,n,\infty]$.  The interior angles of
$T[m,n,\infty]$ are once $\pi/2$ and $m/2$ times $\pi/n$, in this
order.  The remaining angle is $(3n-m)\pi/2n$ if $m\equiv 0\bmod{4}$
and $((n-m)\pi/(2n)$ if $m\equiv 2\bmod{4}$. We conclude that the
number of self-crossings is $(m-4)/4$ (resp.\ $(m-2)/4$) if $m\equiv
0\bmod{4}$ (resp.\ $m\equiv 2\bmod{4}$).  Therefore the number of
self-crossings is zero if and only if $m=2,4$.  The case $m=2$
corresponds to Veech's family \cite{Ve89} (\S \ref{ngonrev}). We
show below that the case $m=4$ corresponds to the billiards constructed
in \S \ref{4nsec}.

 We leave it to the reader to use Theorem \ref{hisTeich} and the
techniques of Theorem \ref{Wardexplicit} to construct billiard tables
with projective affine group $\Delta(4,n,\infty)$ and
$\Delta(5,n,\infty)$ also in the case that  $n$ even or  divisible by
$5$, or both.  

\begin{Prop}\label{Betaprop}
Let $m$ be either $4$ or $5$.
 The billiard table $T[m,n,\infty]$ is similar to the
billiard table $T{(m,n,\infty)}$.
\end{Prop}
\par
{\bf Proof:} 
Suppose that $m=5$. The case that $m=4$ is similar, and left to the reader.

Recall that the interior angles of the $4$-gons $T(5,n,\infty)$ and
$T[5,n, \infty]$ are the same, and also occur in the same
order. Therefore we only have to compare the lengths of the sides of
$T[5,n,\infty]$ to those of $T(5,m,\infty)$. Since the sides of $T[5,
n, \infty]$ are expressed in terms of the Schwarz--Christoffel map, it
suffices to show that
\begin{equation} \label{SCidentity}
 \frac{|SC(2\cos(2\pi/5))-SC(2)|}{|SC(\infty)-SC(2)|} =
\frac{|I_3|}{|I_4|} = \frac{\cos(\pi/n)+\cos(\pi/5)}{\cos(\pi/2n)}.
\end{equation}
Here $I_3, I_4$ are the vectors corresponding to the sides of the
$4$-gon $T[5,n,\infty]$ as indicated in Figure \ref{billiard59}.

We first express the length of the vector $I_4$ in terms of Beta integrals:
\begin{equation}
\begin{split}
|I_4|&=\int_2^\infty (u-2)^{\frac{1}{2n}-1}(u-2\cos(2\pi/5))^{\frac{1}{n}-1}
(u-2\cos(4\pi/5))^{\frac{1}{n}-1}\, {\rm d} u\\
&= \int_1^\infty z^{1-\frac{5}{2n}}(z^5-1)^{\frac{1}{n}-1}(z+1)\, {\rm d} z.
\end{split}
\end{equation}
Here we used the substitution $u=z+1/z$, compare to the proof of
Theorem \ref{Wardexplicit}. Substituting $z=1/t$, we recognize this
integral as the sum of two Beta integrals:
\begin{equation}\label{I4eq}
|I_4|=\frac{1}{5}\left(\mathop{\rm B}(\frac{2}{5}-\frac{1}{2n}, \frac{1}{n})+
\mathop{\rm B}(\frac{3}{5}-\frac{1}{2n}, \frac{1}{n})\right).
\end{equation}

Similarly, one finds that
\begin{equation}\label{I3eq}
|I_3|=\frac{1}{5}\left((-1+\zeta_5^2\zeta_{2n}^{-1}) \mathop{\rm
 B}(\frac{2}{5}-\frac{1}{2n},
 \frac{1}{n})+(-1+\zeta_5^3\zeta_{2n}^{-1}) \mathop{\rm
 B}(\frac{3}{5}-\frac{1}{2n}, \frac{1}{n})\right).
\end{equation}
Equation (\ref{SCidentity}) follows from (\ref{I4eq}) and (\ref{I3eq})
by expressing the Beta integrals in terms of Gamma functions, using
that $\Gamma(z)\Gamma(1-z)=\pi/\sin(\pi z)$, and applying the addition
formulas for sines and cosines.
\hspace*{\fill} $\Box$

One may give an alternative proof for the statement that the $4$-gons
$T[5,n,\infty]$ and $T(5,n,\infty)$ are similar by showing the following.  Let
$P$ be any $4$-gon with the prescribed angles, and let $X$ be the
corresponding translation surface. Suppose that the affine group of $X$
contains $R$ and $T$. Then $P$ is similar to $T(5,n,\infty)$. This can be
shown by first deducing from the geometry of $P$ that an affine diffeomorphism
with derivative $R$ has to fix each saddle connection.

\par
\begin{Rem} {\rm
Several authors (\cite{Wa98}, \cite{Vo96}, \cite{KeSm00}, \cite{Pu01})
have classified the Teich\-m\"uller curves that are obtained by
unfolding a rational triangle, under certain conditions on the angles
of the triangle. We have obtained the translation surfaces $X(m,n,\infty)$
for $m=4,5$ by unfolding $4$-gons. The corresponding families of
Teichm\"uller curves have not been found by Ward et.\ al.\ This
suggests that the translation surfaces $X(m,n,\infty)$ for $m=4,5$
may not be obtained by unfolding triangles, but of course we have not
shown this.  
}
\end{Rem}

\begin{Rem}{\rm
For $n>m\geq 6$ we have not been able to obtain the translation
surface $X(m,n,\infty)$ by unfolding a billiard table, since the
corresponding polygon $T[m,n,\infty]$ may not be embedded in the
complex plane. However, it should in principle be possible to give a
concrete description of $X(m,n,\infty)$ as obtained by gluing certain
cylinders, analoguous to the description in the case of $m=4,5$
(\S\S \ref{5nsec}, \ref{4nsec}). As for $m=4$ and $5$, it follows from
Corollary \ref{splitVHSh} that we would need $g(X_0)$ cylinders, which
is approximately $(m-1)(n-1)/2$: it will be difficult to
visualize the result. Therefore it seems more natural to us to
represent these Teichm\"uller curves via the algebraic description from
\S \ref{realize}.
}
\end{Rem}


\section{Lyapunov exponents} \label{Lyapunov}

Roughly speaking, a flat normed vector bundle on a manifold with
a flow, i.e.\ an action of $\RR^+$, can sometimes be
stratified according to the growth rate of the length of
vectors under parallel transport along the flow.
The growth rates are then called Lyapunov exponents. In this
section we will relate Lyapunov exponents to degrees of
some line bundles in case that the underlying manifold is a
Teichm\"uller curve.
\newline
For the convenience of the reader we reproduce Oseledec's
theorem (\cite{Os68}) that proves the existence of such exponents.
We give a restatement due to \cite{Ko97} in a language closer
to our setting. 
\par

\subsection{Multiplicative ergodic theorem}

We start with some definitions.
A measurable vector bundle is a bundle that can be trivialized
by functions which only need to be measurable. If $(V,||\cdot||)$
and $(V',||\cdot||)$ are a normed vector bundles and 
$T: V \to V'$, then we let $||T||:= \sup_{||v||=1} ||T(v)||$. 
A reference for notions in ergodic theory is \cite{CFS82}.
\par
\begin{Thm}[Oseledec]
Let $T_t: (M,\nu) \to (M,\nu)$ be an ergodic
flow on a space $M$ with finite measure $\nu$.
Suppose that the action of $t \in \RR^+$ lifts
equivariantly to a flow $S_t$ on some measurable 
real bundle $V$ on $M$. Suppose there exists
a (not equivariant) norm
$||\cdot||$ on $V$ such that for all $t\in\RR^+$
$$ \int_M \log(1+||S_t||)\nu < \infty.$$ 
Then there exist real constants $\lambda_1 \geq \cdots \geq \lambda_k$
and a filtration 
$$V = V_{\lambda_1} \supset \cdots V_{\lambda_k} \supset 0$$
by measurable vector subbundles such that, for almost all $m \in M$ 
and all $v \in V_m \sms \{0\}$, one has
$$||S_t(v)|| = {\rm exp}(\lambda_i t + o(t)),$$
where $i$ is the maximal value such that
$v\in (V_i)_m$.
\newline
The $V_{\lambda_i}$ do not change if $||\cdot||$
is replaced by another norm of `comparable' size (e.g.\
if one is a scalar multiple of the other).
\end{Thm}
\par
The numbers $\lambda_i$ for $i=1,\ldots,k\leq {\rm rank}(V)$
are called  the {\em  Lyapunov exponents of $S_t$}. Note
that these exponents are unchanged if we replace
$M$ by a finite unramified covering with
a lift of the flow and the pullback of $V$. We adopt the convention
to repeat the exponents according to the
rank of $V_i/V_{i+1}$ such that we will always
have $2g$ of them, possibly some of them equal. 
A reference for elementary properties of
Lyapunov exponents is e.g.\ \cite{BGGS80}.
\par
If the bundle $V$ comes with a symplectic structure 
the Lyapunov exponents are symmetric with respect to
$0$, i.e.\ they are (\cite{BGGS80} Prop.\ 5.1)
$$ 1= \lambda_1 \geq \lambda_2\geq \cdots \geq \lambda_g \geq 0
\geq -\lambda_g \geq \cdots \geq -\lambda_1=-1.$$
\par
We specialize these concepts to the situation we are
interested in. Let $\Omega M_g^*$ be the bundle of
non-zero holomorphic $1$-forms over the moduli space of curves.
Its points are translation surfaces. The $1$-forms define
a flat metric on the underlying Riemann surface
and we let $\Omega_1 M_g \subset \Omega M_g^*$ 
be the hypersurface consisting of translation surfaces of area one.
As usual we replace $M_g$ by an appropriate fine moduli space
adding a level structure, but we do not indicate this in the
notation. This allows us to use a universal
family $f: \cXX \to M_g$.
\newline
Over $\Omega_1 M_g$, we have the local system 
$\VV_\RR = R^1 f_* \RR$, whose fiber over $(X,\omega)$ is 
$H^1(X,\RR)$. We denote the corresponding real $C^{\infty}$-bundle by $V$.
This bundle naturally carries the Hodge metric
$$H(\alpha,\beta) = \int_X \alpha\wedge*\beta,$$
where classes in $H^1(X,\RR)$ are represented by $\RR$-valued
$1$-forms, and where $*$ is the Hodge star operator. We denote
by $||\cdot||:= ||\cdot||_T$ the associated metric on $V$.
\newline
There is a natural $\SL_2(\RR)$-action on $\Omega_1 M_g$ 
obtained by post-composing the charts given by integrating
the $1$-form with the $\RR$-linear map given by $A \in \SL_2(\RR)$
to obtain a new complex
structure and new holomorphic $1$-form (see e.g.\  \cite{McM03}
and the reference there). 
The geodesic flow $T_t$ on $\Omega_1 M_g$ is the restriction
of the $\SL_2(\RR)$-action to the subgroup ${\rm diag}(e^t, e^{-t})$. 
Since $V$ carries a flat structure, we can lift $T_t$ by parallel 
transport to a flow $S_t$ on $V$. This is the {\em Kontsevich--Zorich cocycle}.
The notion `cocycle' is motivated by writing the flow on a vector bundle
in terms of transition matrices.
\par
Lyapunov exponents can be studied for any finite measure $\nu$ on a 
subspace $M$ of $\Omega_1 M_g$ such that $T_t$ is 
ergodic with respect to $\nu$. Starting with the work
of Zorich (\cite{Zo96}), Lyapunov exponents have been studied for
connected components of the stratification of $\Omega_1 M_g$
by the order of zeros of the $1$-form. The integral structure 
of $\Omega M_g^*$ as an affine manifold can be used to construct 
a finite ergodic measure $\mu$.
Lyapunov exponents for $(\Omega_1 M_g,\mu)$ may be interpreted as deviations
from ergodic averages of typical leaves of measured foliations
on surfaces of genus $g$. The reader is referred to 
\cite{Ko97}, \cite{Fo02} and the surveys \cite{Kr03} and \cite{Fo05}
for further motivation and results.
\par
\subsection{Lyapunov exponents for Teichm\"uller curves}

We want to study Lyapunov exponents in case of 
an arbitrary Teichm\"uller curve $C$ or rather its canonical lift $M$
to $\Omega_1 M_g$ given by providing the Riemann surfaces 
parameterized by $C$ with the normalized generating differential. 
The lift $\pi: M \to C$ is an $S^1$-bundle.
We equip $M$ with the measure $\nu$ which is induced by
the Haar measure on $\SL_2(\RR)$, normalized such that $\nu(M)=1$.
Locally, $\nu$ is the product of the measure $\nu_C$ coming 
from the Poincar\'e volume form and the uniform measure on $S^1$, 
both normalized to have total volume one.
\newline
We can apply Oseledec's theorem since $\nu_M$ is ergodic
for the geodesic flow (\cite{CFS82} Theorem~4.2.1).
\par
We start from the observation that the decomposition $(\ref{VHSdecomp})$
of the VHS in Theorem\ \ref{maxHiggsTeich}
is $\SL_2(\RR)$-equivariant and orthogonal with respect to Hodge metric.
This implies that the Lyapunov exponents of $\VV$ are the
union of the Lyapunov exponents of the $\LL_i$ with those of $\MM$.
\par
Let $\cLL_i:= (\LL_i)^{1,0}$ be the $(1,0)$-part of the
Hodge filtration of the Deligne extension of $\LL_i$ to $\ol{C}$.
Denote by $d_i := \deg(\cLL_i)$ the corresponding degrees.
Recall from Theorem \ref{maxHiggsTeich} that precisely one of the $\LL_i$, say 
the first one $\LL_1$ is maximal Higgs. Recall that $S = \ol{C} \sms C$
is the set of singular fibers.
\par
\begin{Thm} \label{rLyap}
Let $\nu_M $ be the finite $\SL_2(\RR)$-invariant 
measure with support in the canonical lift $M$
of a Teichm\"uller curve to $\Omega_1 M_g$. Then
$r$ of the Lyapunov exponents $\lambda_i$ satisfy 
$$\lambda_i = d_i/d_1 = \lambda(\LL_i,S).$$
In particular, these exponents are rational, non-zero
and their denominator is bounded by $2g-2+s$, where $s = |S|$. 
\end{Thm}
\par
{\bf Proof:} We write $(\LL_i)_\RR$ for the local subsystem of $R^1 f_* \RR$
such that $(\LL_i)_\RR \otimes_\RR \CC = \LL_i$ and let $L_i$
be the $C^{\infty}$-bundle attached to $(\LL_i)_\RR$.
We apply Oseledec's theorem to $L_i$. Then 
$$\lambda_i = \lim_{t \to \infty} \frac{1}{t} \log||S_t(v_i)||,$$
for $v_i \in L_i \sms (L_i)_{-\lambda_i}$.  By averaging, we have
$$\lambda_i = \lim_{t \to \infty} \frac{1}{t} \int_{G(L_i)}
\log||S_t(v_i)|| \dd\nu_{G(L_i)(v_i)},$$
where $\tau: G(L_i) \to M$ is the (Grassmann) bundle of norm one vectors in $L_i$. 
This bundle is locally isomorphic to $S^1 \times M$. The measure  $\nu_{G(L_i)}$ is
locally the product measure of $\nu$ with the uniform measure on $S^1$.
\par
Following the idea of Kontsevich (\cite{Ko97})
also exploited in Forni (\cite{Fo02}), we 
estimate the growth of the length of $v_i$ not only as a
function on the $T_t$-ray through $\tau(v_i)$ (given as
the parallel transport of the corresponding vector) but as
a function on the whole (quotient by a discrete group of a) 
Poincar\'e disc $D_{\tau(v_i)}$ in  $M$. For this purpose
we write $z= e^{i\theta} r $ ($\theta \in [0,2\pi))$ for $z$ in the unit disc 
$D$ and lift it to $\rho_\theta \, \diag(e^t, e^{-t}) \in \SL_2(\RR)$, 
where $\rho_\theta$ is the rotation matrix by $\Theta$ and 
$t =(1/2)\log((1+r)/(1-r))$. 
Using this lift $D \to \SL_2(\RR)$ we obtain our disc
$D_{\tau(v_i)}$ in  $M$ using the (left) $\SL_2(\RR)$-action
on $M$.
\newline
Consider the following functions  
$$ f_D := f_{D,i}:  \left\{ \begin{array}{lcl} 
(\pi^* L_i \sms \{0\}) \times D  & \longrightarrow &   \RR \\ 
(v_i,z) & \mapsto & \log ||z \cdot v_i ||, \end{array} \right. $$
where  $z \cdot v_i$ is the parallel transport of $v_i$ over the 
disc $D_{\tau(v_i)}$.  This is well-defined since
the monodromy of $L_i$ acts by matrices in $\SL_2(\ZZ) = {\rm Sp}_2(\ZZ)$
and symplectic transformations do not affect the Hodge length. 
Note that by definition 
\begin{equation} \label{ftrans}
f_D(v_i,z) = f_D(z \cdot v_i,0). 
\end{equation}
\par
On the discs $D_{\tau(v_i)}$ we may apply the (hyperbolic) Laplacian $\Delta_h$
to the functions $f_{D_{\tau(v_i)}}$ with respect to the second
variable, i.e.\ consider
$$ h_D := h_{D,i} : \left\{ \begin{array}{lcl} 
(\pi^* L_i \sms \{0\}) \times D & \longrightarrow & 
\RR \\ (v_i, z) & \mapsto & (\Delta_h f_D(v_i,\cdot))(z).
\end{array} \right. $$
Using (\ref{ftrans}) and the invariance of $\Delta_h$ under isometries
one deduces that there is a function $h: \pi^* L_i \sms \{0\} \to \RR$, such
that 
\begin{equation}
h_D(v_i,z) = h(z \cdot v_i).
\end{equation}
\par
Since obviously $\int_{G(L_i)} h(S_t v_i) \dd \nu_{G(L_i)}(v_i) = 
\int_{G(L_i)} h(v_i) \dd \nu_{G(L_i)}(v_i)$ for any $t$, we can apply 
\cite{Kr03} Equation (3) (see also \cite{Fo02}
Lemma 3.1) to obtain
\begin{equation} \label{lyviah}
\lambda_i = \int_{G(L_i)}  h(v_i) \nu_{G(L_i)}(v_i). 
\end{equation}
\par
We want to relate this expression to the degree $d_i$ of the
line bundles $\cLL_i$. Suppose $s_i(u)$ is a holomorphic section
of $\cLL_i$ over some open $U \subset C$. Recall that $L_i$ has 
unipotent monodromies, by assumption. Therefore \cite{Pe84} Proposition\ 3.4
implies that the Hodge metric grows not too
fast near the punctures and we have
\begin{equation} \label{diviacurvature}
d_i = \frac{1}{2\pi i} \int_{\ol{C}} \partial \ol{\partial} 
\log(|| s_i||).
\end{equation}
Here as usual, if there is no global section of $\cLL_i$ the contributions of 
local holomorphic sections are added up using a partition of unity.
\par
Instead of considering a holomorphic section $s_i$, we now consider
a flat section $v_i(u)$ of $L_i$ over $U$. Then, in 
$(\wedge^2(\LL_i)_\CC)^{\otimes 2}(U)$
one checks the identity
\begin{equation}
(v_i\wedge*v_i)\otimes(s_i\wedge \ol{s_i}) = \frac{1}{2} 
(v_i \wedge s_i)\otimes
(v_i \wedge \ol{s_i}). 
\end{equation}
\par
We integrate this identity over the fibers $\cXX_c$ of $f: \cXX \to C$, take
logarithms and the Laplacian $\frac{1}{2\pi i} \partial
\ol{\partial}$. Note that
\begin{equation} \label{letzteformel}
\frac{1}{2\pi i} \partial \ol{\partial} \log \frac{1}{2} 
(v_i \wedge s_i)\otimes
(v_i \wedge \ol{s_i}) = 0.
\end{equation}
Let $F$ be a fundamental domain for the action of the affine group $\Gamma$ 
in a Poincar\'e discs $D \hookrightarrow M$. Then (\ref{diviacurvature}) 
and (\ref{letzteformel})
implies that for any flat section $v_i$ of $L_i$ we have
$$  
d_i = \frac{-1}{2 \pi i} \int_F \partial \ol{\partial} 
\log(||v_i||).
$$
The differential operator $\partial \ol{\partial}$ 
coincides, up to a scalar, with $\Delta_h(\cdot) \omega_P$, 
where $\omega_P$ is the Poincar\'e area form. Therefore we obtain
for each $v_i \in (\pi\circ \tau)^* (L_i \sms \{0\})$ that
$$
d_i = \frac{1}{4\pi}\int_F \Delta_h \log ||v_i(z) || \omega_P(z),
$$
where $v_i(z)$ is obtained from $v_i$ via parallel transport.
Hence by integrating over all $G(L_i)$ and taking care of
the normalization of $\nu_{G(L_i)}$ we find that
\begin{equation}\label{bla}
d_i = \frac{1}{4\pi}{\rm vol}(C) \int_{G(L_i)} \Delta_h \log ||v_i|| 
\nu_{G(L_i)(v_i)}
\end{equation}
The statement of the theorem now follows by comparing 
(\ref{bla}) with (\ref{lyviah}). \hspace*{\fill} $\Box$
\par
\begin{Cor} At least $r$ of the Lyapunov exponents are non-zero.
\end{Cor}
\par
{\bf Proof:} By Theorem \ref{rLyap}, it is sufficient to show that
for $\cLL_i := (\LL_i)^{(1,0)}$ the degree $\deg(\cLL_i) \neq 0$.
If $\cLL_i =0$ then, by Simpson's correspondence (\cite{ViZu04}
Theorem 1.1), $\LL_i$ would be a reducible local system. But
since $\LL_i$ is Galois conjugate to $\LL_1$, this is a contradiction. 
\hspace*{\fill} $\Box$.
\par
\begin{Rem}{\rm If $r \geq g-1$ all the Lyapunov exponents are
known. In fact in this case we can identify the remaining
Lyapunov exponent by the formula (\cite{Ko97}, \cite{Fo02} Lemma 5.3)
$$ \sum_{i=1}^g \lambda_i = \frac{\deg(f_* \omega_{X/C})}{2g-2+s} $$
}\end{Rem}
\par
In the case of Teichm\"uller curves associated with triangle
groups constructed in \S \ref{ngonrev} and \S \ref{realize}, 
the proof of Theorem \ref{rLyap} 
yields more. Since for these curves the VHS decomposes
completely into subsystems of rank two (Remark \ref{splitcompletely})
we can determine all the Lyapunov exponents.
\par
\begin{Prop} \label{otherLyap}
Suppose the local system $\MM$ as in Theorem \ref{maxHiggsTeich}
contains a rank two local subsystem $\FF_i$, whose $(1,0)$-part
is a line bundle, which denote by $\cFF_i$. Then the Lyapunov
spectrum contains (in addition to the $d_i/d_1$) the exponents
$$ \deg(\cFF_i)/d_1.$$  
\end{Prop}
\par
By Theorem\ \ref{rLyap} and Proposition\ \ref{otherLyap} it is 
justified to call $\lambda(\LL_i)$ Lyapunov exponents.
\par


\par
Irene I.~Bouw \hfill  Martin M{\"o}ller \newline 
Mathematisches Institut \hfill FB 6 (Mathematik), Campus Essen \newline 
Heinrich-Heine-Universit{\"a}t \hfill Universit{\"a}t Duisburg-Essen \newline
40225 D\"usseldorf \hfill 45117 Essen \newline
bouw@math.uni-duesseldorf.de \hfill martin.moeller@uni-essen.de \newline

\end{document}